\documentclass[12pt,a4paper]{article}

\pdfoutput=1
\usepackage{pdflscape}
\usepackage{amsmath}
\usepackage{amssymb}
\usepackage{latexsym}
\usepackage{srcltx}
\usepackage{graphics}
\usepackage{color}
\usepackage{epsfig}
\usepackage{color}

\newcommand{\re}{{\mathbb R}}
\newcommand{\n}{{\mathbb N}}

\newcommand{\cA}{{\mathcal{A}}}

\newcommand{\cI}{{\mathcal{I}}}

\newcommand{\cL}{{\mathcal{L}}}

\newcommand{\cH}{{\mathcal{H}}}
\newcommand{\cS}{{\mathcal{S}}}

\newcommand{\bx}{{\boldsymbol{x}}}
\newcommand{\by}{{\boldsymbol{y}}}
\newcommand{\bz}{{\boldsymbol{z}}}
\newcommand{\be}{{\boldsymbol{e}}}
\newcommand{\bp}{{\boldsymbol{p}}}

\newcommand{\ba}{{\boldsymbol{a}}}
\newcommand{\bb}{{\boldsymbol{b}}}
\newcommand{\bc}{{\boldsymbol{c}}}
\newcommand{\bu}{{\boldsymbol{u}}}
\newcommand{\bv}{{\boldsymbol{v}}}
\newcommand{\bw}{{\boldsymbol{w}}}
\newcommand{\bm}{{\boldsymbol{m}}}
\newcommand{\bn}{{\boldsymbol{n}}}
\newcommand{\bE}{{\boldsymbol{E}}}

\newtheorem{theorem}{Theorem}
\newtheorem{prop}{Proposition}
\newtheorem{lemma}{Lemma}
\newtheorem{cor}{Corollary}
\newtheorem{remark}{Remark}
\newtheorem{ex}{Example}
\newtheorem{defi}{Definition}

\date{}

\author{
Vladimir Yu.~Protasov
\thanks{Steklov Mathematical Institute of Russian Academy of Sciences, Moscow, Russia,   {e-mail: \tt\small
v-protassov@yandex.ru}}}

\title{Antinorms on cones: duality and applications
\thanks{
This work is supported by the Russian Science Foundation under grant 20-11-20169
}}

\begin{document}
\maketitle 

\begin{abstract}

An antinorm is a concave nonnegative homogeneous functional on a convex cone. 
It is shown that if the cone is polyhedral, then every antinorm  has a unique continuous extension from the interior of the cone. The main facts of the duality theory in convex analysis, in particular,  the Fenchel\,-\,Moreau theorem,  
are generalized  to antinorms.  However, it is shown that the duality relation for antinorms 
is discontinuous. In every dimension there are infinitely many self-dual antinorms on the positive orthant and, in particular, infinitely many 
 autopolar polyhedra. For the two-dimensional case, we characterise them all. 
The classification in higher dimensions is left as an open problem.  Applications to linear dynamical systems, to the Lyapunov exponent of random matrix products,
 to the lower spectral radius of nonnegative matrices, and to convex trigonometry are considered.

\bigskip

\smallskip

\begin{flushright}
\noindent  \textbf{AMS 2000} {\em Mathematical Subject
classification:  46B20, 52A21, 93D20, 37H15}
\end{flushright}

\end{abstract}

\begin{center}
\large{\textbf{1. Introduction}}
\end{center}
\bigskip 

Let $K$ be a convex  cone in~$\re^d$. 
We always assume that a cone is closed, nondegenerate, i.e. possesses a nonempty interior, and
pointed, i.e., does not contain a straight line. Every ray in the cone starting at the apex
is either a  {\em generatrix}, if it lies on the boundary, or an {\em interior ray}
otherwise.

\begin{defi}\label{d.10} 
An {\em antinorm} on a cone $K$ is a nonnegative, somewhere positive, concave
homogeneous functional on $K$. An antinorm is called positive if it is strictly positive at all
points $x \in K\setminus \{0\}$.
\end{defi}
In most cases we assume that $K$ is the positive orthant~$\re^d_+$, but all definitions will be given 
in the general case. 

Antinorms are  concave analogues of norms. However, they cannot be defined in the whole space because there are no
concave positively homogeneous functions on $\re^d$. 
That is why, an antinorm is  usually restricted to some   cone. 
It follows from the concavity that an antinorm is strictly positive (i.e., does not vanish) in the interior of $K$. 
It can vanish at some points of the boundary  $\partial K$, which makes a difference form norms. 
Another difference is that 
antinorms can be discontinuous. However, as we shall see in Section~2, for polyhedral cones, in particular, for~$\re^d_+$, 
the whole theory can be reduced to continuous antinorms by using the concept 
of {\em continuous extension}.  The most remarkable differences between norms and antinorms are in the duality theory, which is  
 developed  in Sections~3.  
 The main concepts and the basic facts (duals and polars, Young's inequality, refrexivity of duality, etc.) are similar to those for dual norms.   
 However, the duality relation can be  discontinuous even for 
 continuous antinorms (Theorem~\ref{th.15} in Section~3). Moreover,
 while  a self-dual norm in~$\re^d$ is unique (the Euclidean one), there exist infinitely 
 many different self-dual antinorms in~$\re^d_+$. In case~$d=2$, we give their 
 complete classification. On the other hand, there is only one symmetric self-dual 
 antinorm in~$\re^2_+$, which is bounded by a hyperbola. This issue is studied in Section~4.  
 Generalizations of those results to higher dimensions are left as open problems in Section~5.
 In Section~6   we analyse applications of antinorms to positive linear switching systems, to asymptotics of random matrix products, to the lower spectral radius, and to convex trigonometry. 
 \smallskip

  \textbf{Related works}. Generalizations of the notion of the norm such as 
  pseudonorm, seminorm, etc. have been thoroughly studied in the literature. 
  The concept of antinorm is usually understood according to 
  Definition~\ref{d.10}, as a positively homogeneous functional 
  (i.e., $f \ge 0$ and $f(\lambda \bx) = \lambda f(\bx)$ for all 
$\lambda \ge 0$)
  with the reverse triangle inequality. As a rule, an antinorm is  defined on a cone or, more 
  generally, is defined piecewise on a 
  fan, which is a partition of a space into several cones with a 
  common apex. Similarly antinorms were defined~\cite{MR12, R15}
  as  piecewise concave Minkowski functionals of star sets. 
  This definition was put to  good use for special construction of probability distributions~\cite{R15}. Matrix antinorms  on the cone of positive definite symmetric matrices~$\mathbb{M}_d$ were 
studied in~\cite{BH11,  BH14},  where some important ineqialities for operator means were extended to antinorms. Again, the antinorm was defined as a positively homogeneous 
 concave function (of a matrix).  In particular, 
 the Minkowski antinorm $f(A) = ({\rm det}\, A)^{1/d}$ and 
 the Schatten $q$-antinorms~$f(A) = ({\rm tr}\, A)^{1/q}, \, q \in (-\infty, 1]$, 
 on the set~$\mathbb{M}_d$ were considered in those works. 
  The concept of antinorm was extended to von Neumann algebras~\cite{BH15}. 
 On the other hands, in some works the term ``antinorm'' has a different meaning. 
 For example,~\cite{MS06} a dual norm is called antinorm. In that work 
 dual norms were analysed in context of the Radon curves and of the Minkowski content

 We use the notion of antinorm according to Definition~\ref{d.10} as 
 it is done in most of the literature. 
 To the best of our knowledge, this concept 
 originated with Merikoski~\cite{M90}
 and was studied in \cite{MS91, MO92}. 
 Independently and later  antinorms were defined in aforementioned works. 
  It was shown~\cite{P10} that every i.i.d. sequence of random non-negative matrices, 
 under some mild assumptions, possesses an invariant antinorm on the positive orthant. 
 This result was applied for the problem of computation of the largest 
 Lyapunov exponent for random matrices~\cite{JP13}. Invariant  antinorms were 
 also exploited 
 for computation of the lower spectral radius of matrix families~\cite{GP13, GZ15, GZ20}. The work~\cite{GZ20}  also analyses many properties of antinorms, including the 
 basic facts of their duality. The antinorms 
 were used
 in the stability theory for positive linear dynamical systems~\cite{GLP17}.  
 
   \smallskip 
 
 \begin{remark}\label{r.50}
{\em Since a positive homogeneous function 
cannot be concave  on the whole space~$\re^d$,  
antinorms  are usually defined on convex cones. For example, they are often considered on  the positive orthant~$\re^d_+$ or, in case of matrix antinorms, in the cone of positive semitefinite 
matrices~\cite{BH11,  BH14}.   
In some works antinorms are defined  on~space~$\re^d$, but in this case 
the whole space is split to several cones (forming a fan) 
and the antinorm is concave  on each cone  separately.  An equivalent definition uses Minkowski functionals of star sets~\cite{MR12, R15}.}
\end{remark}
 
  \textbf{Notation}.  
For an antinorm~$f$ on~$K$, one defines its {\em antiball}~$G = \{\bx \in K \ | \ 
 f(\bx) \ge 1\}$ and an {\em antisphere} $\cS = \{\bx \in K \ | \ 
 f(\bx) = 1\}$. The antiball is convex and unbounded. Moreover, it is 
 a {\em conic body} according to the following definition: 
 \begin{defi}\label{d.15} 
A {\em conic body} is a convex closed subset of a cone~$K$ that 
does not contain the apex and possesses the following property: 
every interior ray of~$K$  intersects that set by a ray. 
A {\em conic polyhedron} is a conic body defined in~$K$ by a system of 
finitely many linear inequalities. 
\end{defi}

Clearly, all conic bodies are unbounded. 
Unlike usual polyhedra, a conic polyhedron can have only one vertex, 
for example~$G = \{x \in \re^d_+\ | \ x_i \ge 1, i = 1, \ldots , d\}$. 

\begin{figure}[ht!]
\begin{minipage}[h]{0.3\linewidth}
\center{\includegraphics[width=1\linewidth]{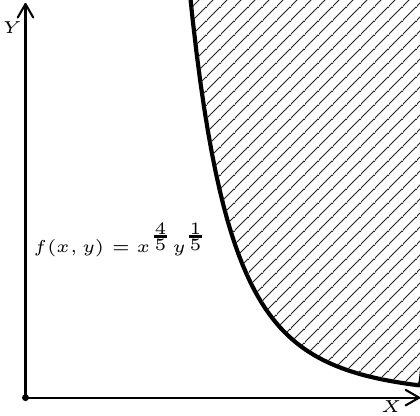}} a) \\
\end{minipage}
\hfill
\begin{minipage}[h]{0.3\linewidth}
\center{\includegraphics[width=1\linewidth]{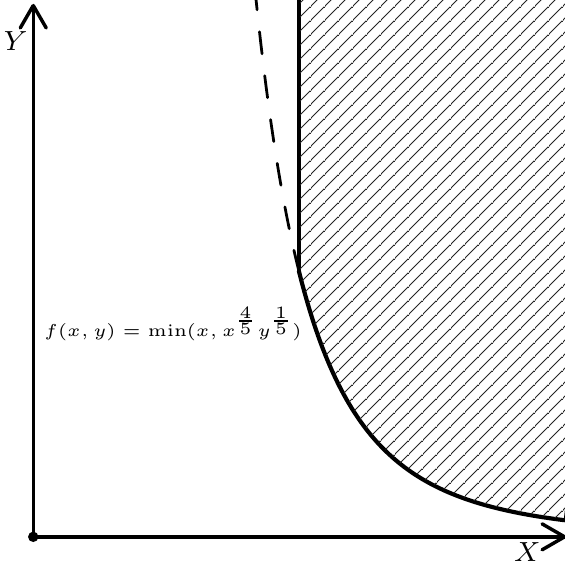}} b) \\
\end{minipage}
\caption{The unit antiballs of the antinorms~$\, x^{\,\frac45}\, y^{\,\frac15}$ (left)
and ~$\min\, \bigl\{ \, x^{\,\frac45}\, y^{\,\frac15} \, , \, x\bigr\}$ (right).}
\label{pic03}
\end{figure}

There is a one-to-one correspondence between 
antinorms and conic bodies. The antiball of an antinorm~$f$ is  a conic body. 
Conversely, every conic body~$G$ defines a unique 
antinorm, similarly to the Minkowski functional: $f(\bx) \, = \, \sup \, \{
\lambda > 0\ | \ \lambda^{-1}\bx \in G \}$. 
The dual antinorm corresponds to the polar conic body. There are infinitely many 
autopolar conic polygons in~$\re^2_+$. In Section~4 we explicitly classify them all. We are not aware 
of any example of autopolar polyhedra in~$\re^d_+$ for $d\ge 3$ except for those 
reduced to two-dimensional ones. 
In what follows we usually drop the prefix~``anti'' is use the simple terminology 
``ball'' in ``sphere'' for antinorms, when it does not lead to confusions. 
Fig.~\ref{pic03} shows the unit balls  of two antinorms.

The vectors are denoted by bold letters and the scalars are denoted by 
usual letters, 
so $\bx = (x_1, \ldots , x_d)$. We use the standard nonation 
$\re^d_+$ for the positive orthant, which  consists of 
points with all coordinates being nonnegative; $\bx \ge 0$ if $\bx \in \re^d_+$
and $\bx > 0$ if $\bx \in {\rm int}\, \re^d_+$; $\bx \ge \by$ means that $\bx - \by \, \ge \, 0$.   

As usual, ${\rm int}\, G, \, \partial \, G$, and ${\rm co}\, G$ denotes respectively the interior, the boundary, and the convex hull  of a set~$G$. A {\em convex body} is a convex compact set with a nonempty interior. 
We denote the Euclidean norm in~$\re^d$ by $|\, \cdot \, |$ and an arbitrary norm 
by~$\|\, \cdot \, \|$. 
By~$B_r(\bx)$ we denote the Euclidean ball of radius $r$ centered at 
 $\bx$. An {\em inverse point} to a given point~$A \in \re^d$ is the image of 
 $A$ under the inversion about the unit sphere centered at the origin. 
 Thus, $A'$ is inverse to $A$ if the vectors $OA$ and $OA'$ are co-directed and 
 $|OA'|\cdot |OA|\, = \, 1$.

\bigskip

\begin{center}
\large{\textbf{2. Continuity of antinorms}}
\end{center}

\medskip 

It is well known that  antinorms may be discontinuous. 
For example, the following antinorm in~$\re^2_+$:  $
 f(x, y) = x+y$ if $x>0$ and $f(x, y) = 0$ if $x=0$,  is discontinuous. 
The continuity issue is extremely important  in 
generalizing many facts of convex analysis to antinorms. For example, 
the main duality result $f^{**} = f$ is true provided an antinorm is continuous and may fail otherwise~\cite{GZ20}. Moreover, continuity of antinorms is crucial in many applications.
For instance, extremal and invariant antinorms of linear dynamical systems must be continuous to  estimate the growth of trajectories, see Section~6.

In case of antinorms on the positive cone~$\re^d_+$, or, more generally, on a polyhedral cone, the whole  theory can be restricted to the continuous case. To see this we need 
 to establish some properties of concave functions defined on polyhedral sets. 

A convex closed set  $G \subset \re^d$ with a nonempty interior is called {\em polyhedral} 
if it is a set of 
solutions of a system of linear inequalities. A bounded polyhedral set is a polyhedron. 

Since a concave function is continuous at any interior point 
of its domain, it follows that the discontinuity 
can occur only on the boundary.  
Let $f$ be a non-negative concave function defined on the interior of a convex set $G$. We extend it onto the boundary of $G$ by the limit:  
\begin{equation}\label{eq.boundary}
F(\bx)\ = \ \lim_{{\by\in {\rm int}\, G \atop |\bx-\by|\to 0}} \, f(\by)\ , \qquad 
\bx \in \partial G\, .
\end{equation}
\begin{lemma}\label{l.10}
If $G$ is a polyhedral set, then for every point $\bx \in \partial G$, the limit~(\ref{eq.boundary}) is well-defined. 
\end{lemma}
{\tt Proof.} Assume the contrary. In this case there are two numbers 
$a, b$ such that $a < b$ and two sequences $\{\by_k\}_{k\in \n}, \, \{\bz_k\}_{k\in \n}$
of points from ${\rm int}\, G$ that converge to $\bx$ and such that 
$f(\by_k) < a < b < f(\bz_k)$ for  all~$k\in \n$. Take a number $r$
so small that the ball~$B_r(\bx)$  does not 
intersect other face planes of $G$ than those containing $\bx$. Denote 
by $M$ the supremum of the function $f$ in the interior of this ball. 
The concavity implies that $M< +\infty$. Let  
$\varepsilon < \frac{r(b-a)}{M}$. For every point $\bz \in 
B_{\varepsilon}(\bx)\cap G$, the point $\bx + \frac{r}{\varepsilon}(\bz-\bx)$ 
belongs to $B_r(\bx)\cap G$. Now take arbitrary point 
$\bz_k \in B_{\varepsilon}(\bx)\cap G$ and a point $\by_n$ very close to $\bx$. 
Then the point $\ba = \bx + \frac{r}{\varepsilon}(\bz_k-\by_n)$ 
belongs to $B_r(\bx)\cap G$. Therefore $f(\ba) \le M$. On the other hand, 
the concavity yields $f(\ba) > f(\bx) + \frac{r}{\varepsilon} (b-a) > 
f(\bx) + M\ge M$. Thus, $f(\ba) > M$, which is a contradiction.  

{\hfill $\Box$}
\medskip 

\begin{remark}\label{r.10}
{\em Lemma~\ref{l.10} may fail for non-polyhedral domains as the following 
example demonstrates.}
\end{remark}

\begin{ex}\label{ex.10}
{\em Consider a disc~$G \subset \re^2$ of radius $1$ centered 
at the point $(1,0)$ and a function~$f$ whose graph in $\re^3$ is a convex 
hull of this disc and of  the point~$(0,0,1)$, see Fig.~\ref{pic01}.  
Take an arbitrary sequence of points $\bx_k$ on the boundary of $G$
which converges to the point $\bx_0 = (0,0)$. Then, obviously, 
both sequences $\frac12\bx_k$ and $\frac13\bx_k$ converge to $\bx_0$. 
However, $\lim_{k\to \infty}f\bigl(\frac12\bx_k\bigr) = \frac12$, while 
$\lim_{k\to \infty}f\bigl(\frac13\bx_k\bigr) = \frac23$. 

\begin{figure}[ht!]
\center{\includegraphics[width=0.5\linewidth]{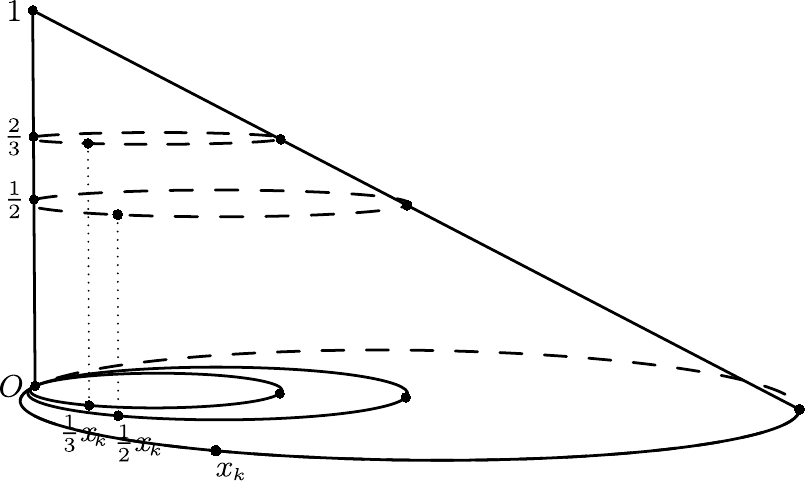}}
\caption{A concave function without continuous extension}
\label{pic01}
\end{figure}

The same construction provides a counterexample to Lemma~\ref{l.10} for 
every non-polyhedral domain. If a domain~$G \subset \re^d$ is not polyhedral, then 
there exist points $\{\bx_k\}_{k=0}^{\infty}$ on its boundary 
such that $\bx_k \to \bx_0$ as $k \to \infty$ and for every $k\ge 1$, 
the open interval between $\bx_0$ and $\bx_k$ lies in the interior of~$G$. 
Then we place $G$ in the coordinate subspace of~$\re^{d+1}$
spanned by the first $d$ basis vectors so that $\bx_0$
coincides with the origin and consider the function~$f$ on $G$ whose graph 
is the convex hull of~$G$ and of the point $(0, \ldots , 0, 1)$. Then we argue as for the case when $G$ is the disc.  
}
\end{ex}

\begin{defi}\label{d.20}
Let $f$ be a nonnegative concave function defined on a convex closed set~$G$
with a nonempty interior. 
A continuous concave function $F$ on~$G$ such that~$F(\bx) = f(\bx)$ for all~$\bx \in 
{\rm int}\, G$ 
is called a continuous extension of~$f$.  
\end{defi}
Not every nonnegative concave function possesses a continuous extension. 
Nevertheless, if the set~$G$ is polyhedral, then such an extension 
does exist. Moreover, it is unique and majorizes  the function~$f$. 
\begin{prop}\label{p.10}
Every nonnegative concave function~$f$ on  a polyhedral set $G$ 
possesses a unique continuous extension~$F$.  That extension is defined on the boundary of $G$ by formula~(\ref{eq.boundary}). For every $\bx\in G$, we have $F(\bx) \, \ge \, f(\bx)$. 
\end{prop}
  
{\tt Proof.} On ${\rm int}\,G$ we have $F=f$. 
At every boundary point the function $F$ is well-defined  by~(\ref{eq.boundary}), as it follows from Lemma~\ref{l.10}. Continuity of~$F$ is easily proved by assuming the 
contrary. Let $F$ be discontinuous at some point~$\bx \in G$. 
Since a concave function is continuous on the interior of 
the domain, it follows that~$\bx \in \partial G$, hence $F(\bx)$
is defined by formula~(\ref{eq.boundary}). There exists a sequence 
$\{\bx_k\}_{k \in \n}$ of points from~$G$, which converges to~$\bx$
but $f(\bx_k)$ does not converge to~$F(\bx)$. Define a new sequence 
$\{\by_k\}_{k \in \n}$ as follows: if $\bx_k \in {\rm int}\,G$, then 
$\by_k = \bx_k$, otherwise $\by_k$ is a point from ${\rm int}\,G$
close to~$\bx_k$ for which $f(\by_k)$ is close to  $F(\bx_k)$
(such points exist due to formula~(\ref{eq.boundary})). 
Then $\by_k \to \bx$, but  $f(\by_k)$ does not converge  to  $F(\bx)$, 
which contradicts to the definition of~$F(\bx)$. 
 
To prove that $F(\bx) \ge f(\bx)$, 
we take an arbitrary point $\by \in {\rm int}\,G$, in which $F(\by) = f(\by)$. 
Concavity implies that 
$$
 (1-t)f(\bx) \ \le \ f\bigl( (1-t)\bx + t\by\bigr)\  - \  t\,f\bigl(\by\bigr) \ = \ 
 F\bigl( (1-t)\bx + t\by\bigr)\  - \  t\,F\bigl(\by\bigr)\, . 
$$
Taking  limit as $t\to +0$, we have 
$ F\bigl( (1-t)\bx + t\by\bigr) \to F(\bx)$ by continuity, and hence $f(\bx) \le F(\bx)$. 
\smallskip

{\hfill $\Box$}

\medskip

If $f$ is continuous, then $F\equiv f$. According to Proposition~\ref{p.10}, all 
concave functions on~$G$ are obtained from some continuous concave function 
by reducing its values on the boundary of~$G$. They can be arbitrarily reduced  
at extreme points, and then reduced in all other points of the boundary 
to keep the function concave.

\begin{cor}\label{c.10}
Every antinorm~$f$ on $\re^d_+$ possesses a unique continuous extension~$F$.
At every  boundary point, we have $F(\bx) \ge f(\bx)$. 
\end{cor}

  Thus, we have proved the following classification 
of antinorms on~$\re^d_+$.
  \begin{theorem}\label{th.5}
Every antinorm in $\re^d_+$ is obtained  from a continuous 
antinorm by arbitrary reducing its values on the boundary 
that keeps the function concave and  homogeneous on all 
coordinate subspaces. 
\end{theorem}
 \begin{ex}\label{ex.11}
 {\em The antinorm~$f: \re^3_+ \to \re_+$ defined by the formula 
$$
 f(\bx) \ = \ 
 \left\{
\begin{array}{lll}
(\sqrt{x_1} + \sqrt{x_2} + \sqrt{x_3})^2 & , & x_1x_2x_3 > 0 \, , \\ 
x_1+ x_2 + x_3 & , & x_1x_2x_3  = 0\, , 
\end{array} 
\right. 
$$
is discontinuous. 
\begin{figure}[ht!]
\begin{minipage}[h]{0.34\linewidth}
\center{\includegraphics[width=1\linewidth]{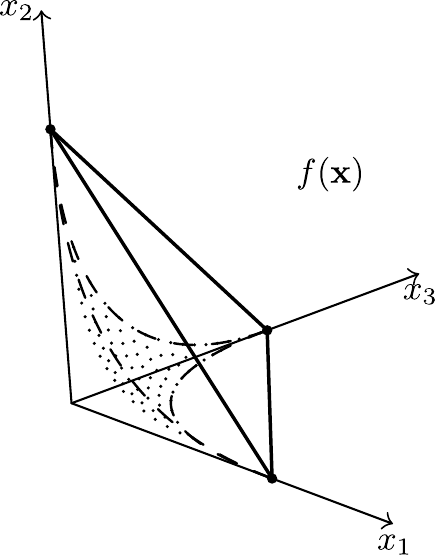}} a) \\
\end{minipage}
\hfill
\begin{minipage}[h]{0.34\linewidth}
\center{\includegraphics[width=1\linewidth]{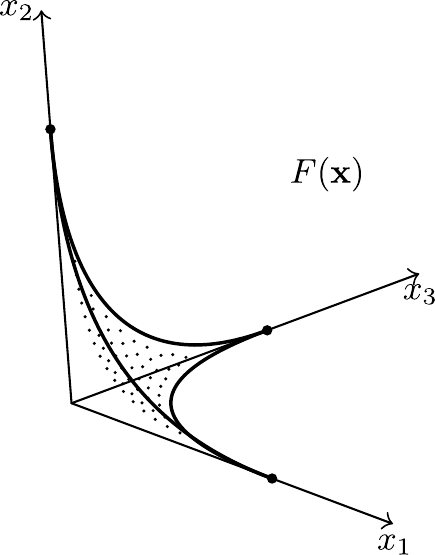}} b) \\
\end{minipage}
\caption{The unit antiball of $f$ (left) and of its continuous extension~$F$ (right) from Example~\ref{ex.11}.}
\label{pic02}
\end{figure}
Its continuous extension is~$F(\bx)\, = \, \bigl( \sqrt{x_1} + \sqrt{x_2} + \sqrt{x_3}
\bigr)^2$. 
On each coordinate subspace, the function~$f$ is obtained from~$F$
by reducing   all values keeping the function concave. At the 
subspace~$\{\bx = (x_1, x_2, 0)\}$, the function~$F = (\sqrt{x_1} + \sqrt{x_2})^2$
is replaced by the smaller concave function $f = x_1 + x_2$, and the same with the other 
coordinate subspaces.    
 }
 \end{ex}

\begin{remark}\label{r.15}
{\em All antinorms on polyhedral cones admit continuous extensions
defined by  formula~(\ref{eq.boundary}). It follows from Proposition~\ref{p.10} that any antinorm on a
polyhedral cone has a continuous extension defined by~(\ref{eq.boundary}). For non-polyhedral
cones
 this may not be true. Indeed, we can consider a Lorentz cone~$K \subset \re^3$, define the 
 antinorm~$f$ on its cross-section disc as in Example~\ref{ex.10} and extend it by homegenity 
 onto the whole~$K$.  Then~$f$ does not have a continuous extension. By a similar 
 argument it can be shown that on every non-polyhedral cone, ther eexists an antinorm without continuous extension.} 
\end{remark}
\smallskip 

In the next section we study duality of antinorms and  prove that an antinorm and its continuous extension have the same dual. This will allow us to replace all antinorms by their extensions,  after which 
we will be able to focus our analysis to continuous antinorms only. 

\bigskip

\begin{center}
\large{\textbf{3. Duality of antinorms}}
\end{center}
\bigskip 

The dual antinorms on arbitrary cones~$K$ were first introduced in~\cite{MS91} and then 
in a more general form in~\cite{GZ15}. 
The duality theory for antinorms is 
very similar to convex duality, with the natural replacement 
of maximum to mininum in the definition of dual functions. 
In particular,  analogues of Young's inequality and of the reflexivity 
for duality of antinorms are true for this case. On the other hand, there are some differences
compared with duality of norms. For example, the double dual 
coincides not with the original antinorm but with its continuous extension. 
Another difference is rather surprising: the correspondence between an antinorm and its dual 
can be discontinuous (even in the class of continuous antinorms). 
\bigskip 

Let us remember that the duality of norms is defined in a standard way: 
if $g(\bx) = \|\bx\|$ is a norm in~$\re^d$, then  $g^*(\bp) = 
\max\limits_{\bx \in \re^d, \bx\ne 0} \frac{(\bp , \bx)}{g(\bx)}$. 
The unit ball $Q^*$ of the dual norm~$g^*$ is a polar to 
the unit ball $Q$ of $g$, i.e., 
$Q^* \, = \, \{\bp \in \re^d \ |\ (\bp , \bx) \le 1, \ \bx \in Q\}$. 

The duals and polars for antinorms are defined in the same way, replacing 
maxima by minima and inverting all inequalities. 
Let $f$ be an antinorm in~$\re^d_+$. Then its dual is 
\begin{equation}\label{eq.dual1}
f^*(\bp)\ = \ \inf_{\bx \in \re^d_+, \, f(\bx) \ne 0} \frac{(\bp , \bx)}{f(\bx)}\ , \qquad 
\bp \in \re^d_+ . 
\end{equation} 
This definition can be modified as follows. If $f(\bx) = (\bp , \bx) = 0$
(this can occur only on the boundary of the cone),
then  $\frac{(\bp , \bx)}{f(\bx)}$ denotes the value  
$\ \liminf_{\by \to \bx} \frac{(\bp , \by)}{f(\by)}$. Being defined this way, the  
function $\frac{(\bp , \bx)}{f(\bx)}$ is lower semicontinuous on~$\re^d_+$,
and hence possesses a point of absolute minimum. This allows us to write 
\begin{equation}\label{eq.dual2}
f^*(\bp)\ = \ \min_{\bx \in \re^d_+} \frac{(\bp , \bx)}{f(\bx)}\ , \qquad 
\bp \in \re^d_+ . 
\end{equation} 
Clearly, the dual function~$f^*(\bp)$ is also an antinorm. The next observation is less trivial: 
an antinorm and its continuous extension  have the same dual. 

\begin{prop}\label{p.15}
If $f$ is an antinorm and $F$ is its continuous extension, then 
$f^* = F^*$. 
\end{prop}
\smallskip 
  
{\tt Proof.} Since $F=f$ on the interior of $\re^d_+$ and 
$F\ge f$ on the boundary, we have 
\begin{equation}\label{eq.chain1}
    \min_{\bx \ge 0} \frac{(\bp , \bx)}{F(\bx)} \ \le \ 
    \min_{\bx \ge 0} \frac{(\bp , \bx)}{f(\bx)} \ \le \ 
    \min_{\bx > 0} \frac{(\bp , \bx)}{f(\bx)} \ = \ 
    \min_{\bx > 0} \frac{(\bp , \bx)}{F(\bx)}\, . 
\end{equation}
On the other hand, the continuity of $F$ implies that 
$\min_{\bx \ge 0} \frac{(\bp , \bx)}{F(\bx)} \, = \, 
\min_{\bx > 0} \frac{(\bp , \bx)}{F(\bx)}$. Hence, all inequalities in 
the chain~(\ref{eq.chain1}) are equalities, therefore $f^*(\bp) = F^*(\bp)$.  
\smallskip 

{\hfill $\Box$}
\medskip

The geometrical meaning of duality is expressed by the {\em antipolar transform}~$X\to X^*$, where 
\begin{equation}\label{eq.antipolar}
X^* \ = \ \Bigl\{ \by \in \re^d \ \Bigl| \  (\by, \bx)\, \ge \, 1 \ \mbox{for all}\ 
\bx \in X \,  \Bigr\} 
\end{equation}
is the {\em antipolar} to the set~$X$. 
We often call it just ``polar transform'' and will avoid confusion with the standard polar transform.  The polar to a point~$\bx$ is the half-space 
$\{\by \in \re^d \ | \ (\bx, \by) \ge 1\}$, and the polar to the set $X$ is the intersection 
of those half-spaces over all~$\bx \in X$. It is shown easily that 
{\em the unit ball of the dual antinorm $f^*$ is the polar to the 
unit ball of~$f$}. 
 \begin{ex}\label{ex.12} 
 {\em For the antinorm~$f(\bx) = \sum_{i=1}^d x_i$, 
 we have 
 $$
 f^*(\bp) \quad = \quad \min_{\bx \in \re^{d}_+, \, \sum_{i=1}^d x_i =1} 
 \ \sum_{i=1}^d p_ix_i\, \quad = \quad \min\, \{p_1, \ldots , p_d\}\, . 
 $$
The unit ball of $f$ is a ``conic simplex'' $G \, = \, \{\bx \in \re^d_+\ | \ \sum_{i=1}^d x_i \ge 1\}$. 
Its polar  $G^* \, = \, \{\bp \in \re^d_+\ | \ p_i \ge 1, \, i = 1, \ldots , d\}$, 
which is the unit ball for $f^*$, is the positive orthant $\re^d_+$
shifted to the vector $\be = (1, \ldots , 1)^T$. 
 }
 \end{ex}

An analogue of Young's inequality for the Legendre-Fenchel transform follows immediately from the definition of dual antinorms. 
\begin{prop}\label{p.20}
For every antinorm and for every points  $\bp, \bx \in \re^d_+$, we have 
\begin{equation}\label{eq.young}
f^*(\bp)f(\bx) \ \ge \ (\bp , \bx)\, .  
\end{equation}
\end{prop}
The  issue of the reflexivity of antinorms was studied in the literature for arbitrary 
convex cones. 
It was shown that~$f^{**} \ge f$~\cite{MS91} . 
Moreover, if an antinorm is continuous, then it is reflexive, i.e.,~${f^{**} = f}$~\cite{GZ20}. It turns out that the latter result can be generalized 
to all antinorms, provided the cone is polyhedral. 
We prove the corresponding theorem for~$\re^d_+$, although it is true for all 
polyhedral cones.

\begin{theorem}\label{th.10}
For an  arbitrary antinorm, we have 
$f^{**} \ = \ F$, where $F$ is the continuous extension of~$F$. 
\end{theorem}
{\tt Proof.} In view of Proposition~\ref{p.15}, it suffices to 
prove that the duality is reflexive for continuous antinorms. 
Indeed, in this case, for an arbitrary antinorm~$f$, we have
$(f^*)^* = (F^*)^* = F^{**} = F$. Therefore, we assume that $f$ is 
continuous and prove that  $f^{**} = f$. To this end we 
need to show that $f(\bx) = \min_{f^*(\bp) \ge 1}
(\bp, \bx) $. Inequality~(\ref{eq.young}) implies that 
$f(\bx) \ge \min_{f^*(\bp) \ge 1}
(p, x)$. It remains to establish the inverse inequality. 
Since $f^*(\bp) = \min_{f(\bz) \ge 1}
(\bp, \bz)$,  the inequality $f^*(\bp) \ge 1$ means that 
 $(\bp, \bz) \ge 1$ for every $z \in B$, where ${B = 
 \{\bz\in \re^d_+ \ | \ f(\bz) \ge 1\}}$.  Thus, 
 $(\bp, \bz) \ge f(\bz)$ for all  $\bz \in \re^d_+$. 
 Hence, we need to show that for every $\bx$ such that 
 $f(\bx) \le 1$ there exists $\bp \in \re^d_+$  such that $(\bp, \bx) = 1$ 
 and $(\bp, \bz) \ge f(\bz)$ for all  $\bz$.
 By the convex separation theorem, there exists $\bp$ such that 
 $1 = (\bp, \bx) \ge \sup_{\bz \in B}(\bp, \bz)$, and hence 
 $(\bp, \bz) \ge f(\bz)$ for all  $\bz$.  It remains to show that 
 $\bp \in \re^d_+$. If this is not the case, $\bp$ is orthogonal 
 to a strictly positive vector~$\ba$. Since  $(\bp, \bz) \, = \, (\bp, \bz + \lambda \ba) \, \ge \, 
 f(\bz + \lambda \ba)$ 
for all $\lambda > 0$. 
 Thus, $\lambda^{-1}(\bp, \bz)  \,  \ge \, 
 f(\lambda^{-1}\bz + \ba)$, which in the limit as $\lambda \to +\infty$
 gives $0 \ge f(\ba)$. This is impossible since $f(\ba) > 0$
 at the positive point~$\ba$. 
 
 {\hfill $\Box$}

\medskip

As a corollary, we obtain the result from~\cite{GZ20}:
\begin{cor}\label{c.20}
If an antinorm~$f$ is continuous, then 
$f^{**} \ = \ f$. 
\end{cor}

Proposition~\ref{p.15} and Theorem~\ref{th.10} allows us to restrict  naturally  our analysis to the set of continuous antinorms. Thus, {\it in what follows, all antinorms are supposed to be continuous if the 
converse is not stated. }
\smallskip 

We see that  the duality map $f\mapsto f^*$ constitutes a reflexive  
transform on the set of antinorms. Surprisingly enough, 
this transform is discontinuous, even on the set of continuous antonorms. This makes one more  difference from 
duality of  norms.

Naturally,  the distance between antimorms is defined as  the 
maximal difference  between them on the unit simplex 
$\Delta = \{\bx \in \re^d_+ \ |  \ 
(\be, \bx) = 1\}$, where $\be$ is the vector of ones. Thus, 
$\|f_1 - f_2\| \, = \, \max_{\bx \in \Delta}|f_1(\bx) - f_2(\bx)| $. 
\smallskip 

\begin{theorem}\label{th.15}
The map $f\mapsto f^*$ is discontinuous. 
\end{theorem}
{\tt Proof}. We present an example in $\re^2$, which is easily extended to higher 
dimensions.  
Consider the following family of antinorms  on~$\re^2_+$: 
$$
f_{\varepsilon}(x, y) \ = \ \min\{x,y\} \ + \ \varepsilon\, \sqrt{xy}\, , \qquad 
 (x,y)\in \re^2_+\,  ,   
$$
where $\, \varepsilon \in [0,1].$
In particular, $f_{0}(x, y) \ = \ \min\{x,y\} $. Let us show that 
$f_{\varepsilon}^*$ does not converge to $f^*_0$ as $\varepsilon \to 0$. 
We have 
$$
f_0^*(p,q) \ = \ \min_{(x,y) \in \re^2_+}\frac{px + qy}{\min\{x,y\}} \ = \ 
\min_{0\le x \le y}\frac{px + qy}{x} \ = \ 
p \, + \, \min_{0\le x \le y}\frac{qy}{x} \ = \ p+q\, . 
$$
Thus, 
\begin{equation}\label{eq.0}
f_0^*(p,q) \ = \ p+q \, .  
\end{equation}
Now compute $f^*_{\varepsilon}(p, q)$ with $\varepsilon > 0$.
We do it only for points $(p,q)$ satisfying the assumption 
\begin{equation}\label{eq.nine}
 q \ \le \ \frac{\varepsilon^2}{8}\, p\, . 
\end{equation}
We have  
$$
f^*_{\varepsilon}(p,q) \ = \ 
\min_{(x,y) \in \re^2_+}\frac{px + qy}{\min\{x,y\} + \varepsilon \sqrt{xy}} \ = \ 
\min_{0\le x \le y}\frac{px + qy}{x  + \varepsilon \sqrt{xy}} \ = \ 
\, \min_{0 < x \le y}\frac{p + q\frac{y}{x}}{1  + \varepsilon \sqrt{\frac{y}{x}}}\, .   
$$
We denote $\sqrt{\frac{y}{x}} = t$ and rewrite as follows: 
$$
f^*_{\varepsilon}(p,q) \ = \ 
\, \min_{t \ge 1}\frac{p + qt^2}{1  + \varepsilon t}. 
$$
Denote this minimum by $m$. This is the smallest parameter value for which 
the equation 
$$
qt^2 \, - \, m\varepsilon t \, + \, p\, - \, m\ = \ 0.   
$$
has a root $t \ge 1$. If $t = 1$, then $m = \frac{p+q}{1+\varepsilon}.$ 
If $t > 1$, then the discriminant
$$
\varepsilon^2 m^2 \, + \, 4 q m \, - \, 4qp \ = \ 0,  
$$  
hence 
$$
m \ = \ \frac{-2q \, + \, \sqrt{4q^2 + 4\varepsilon^2 pq}}{\varepsilon^2} \ = \ 
\frac{4pq}{2q \, + \, \sqrt{4q^2 + 4\varepsilon^2 pq}} \ =  \ 
\frac{2pq}{q \, + \, \sqrt{q^2 + \varepsilon^2 pq}}\, . 
$$
By~(\ref{eq.nine}),  the last expression does not exceed $\frac{p+q}{1+\varepsilon }$. Thus, 
\begin{equation}\label{eq.eps1}
f^*_{\varepsilon}(p,q) \ = \  \frac{2pq}{q \, + \, \sqrt{q^2 + \varepsilon^2 pq}}\, , \qquad q \le \frac{\varepsilon^2}{8}\, p \, .  
\end{equation}
This implies that 
\begin{equation}\label{eq.eps2}
f^*_{\varepsilon}(p,q) \ \le  \  \frac{2}{\varepsilon} \, \sqrt{qp}\, .  
\end{equation}
Therefore, $f^*_{\varepsilon}(1,0) = 0$ for all $\varepsilon$, while 
$f^*_0(1, 0) = 1+ 0 = 1$ as follows from~(\ref{eq.0}). This completes the 
proof. 

{\hfill $\Box$}

\bigskip 

In the next section we characterize self-dualily. We will see that this issue for antinorms is more interesting and challenging 
than for norms. 

\bigskip 

\begin{center}
\large{\textbf{4. Self-dual antinorms}}
\end{center}

\bigskip 

An antinorm is called self-dual if $f^* = f$. Before studying self-dual 
antinorms, let us remember that the situation with self-dual norms is very simple.   
Self-duality is an exclusive property of the Euclidean norm. 
\smallskip 

\noindent \textbf{Fact.} {\em \, The unique self-dual norm 
in $\re^d$ is the Euclidean norm. }
\smallskip 

\noindent For convenience of the reader, we include the proof of this classical fact. 
\smallskip 

 {\tt Proof}. Assume there exists a point 
$\bx \in \re^d$ whose self-dual norm $f(\bx)$ is strictly 
smaller than its Euclidean norm $|\bx|$. Then $f^*(\bx) = \max_{\bz \in \re^d\setminus \{0\}} \frac{(\bx, \bz)}{f(\bz)}\, \ge \,  \frac{(\bx, \bx)}{f(\bx)}\, > \, 
\frac{(\bx, \bx)}{\|\bx\|}\, = \, |\bx|$. Hence, $f(\bx) = f^*(\bx) > |\bx|$, which 
contradicts to the assumption. Thus, $f(\bx) \ge |\bx|$ for all~$\bx$. 
Assume  $f(\bx) > |\bx|$ for some $\bx$. 
In this case $f^*(\bx) > |\bx|$, and hence there exists $\bz \in \re^d_+$
such that $\frac{(\bx, \bz)}{f(\bz)} > |\bx|$. Since ${(\bx, \bz)} \le 
\|\bx\|\cdot \|\bz\|$, it follows that $|\bx|\cdot |\bz| \, > \, |\bx| \cdot f(\bz)$, 
therefore $|\bz| \, > \, f(\bz)$, which is a contradiction.

{\hfill $\Box$}
\medskip 

\begin{remark}\label{r.20}
{\em The statement above means that the unit Euclidean  ball is the only autopolar set 
in~$\re^d$. One may wonder about autopolar triangles or autopolar  simplices, which are well-known. All of them are actually not autopolar: 
their ``autopolarity'' means that the set of vertices is polar to the set 
of sides (or, in~$\re^d$, the set of hyperfaces).   
}
\end{remark}
\medskip 

In contrast, in every dimension~$d \ge 2$,  there are infinitely many different self-dual antinorms, or, which is the same, infinitely many autopolar conic bodies.  
One of families of such antinorms is provided by the following 
assertion. 
\begin{prop}\label{p.30}
Let $\{p_i\}_{i=1}^d$ be a collection of non-negative numbers 
such that   $\sum_{i=1}^d p_i = 1$. Then the function  
$f(\bx) \, = \, \sqrt{d}\, \prod_{i=1}^d\,  x_i^{\,p_i}$ is a self-dual antinorm. 
\end{prop} 
{\tt Proof}. The concavity is well-known. We find $f^*(x)$
 from the problem 
$$
\left\{
\begin{array}{ccc}
- \sqrt{d}\, \prod_{i=1}^d z_i^{\,p_i} & \to & \min\\
\sum_{i = 1}^d x_iz_i & = & 1\, . 
\end{array}
\right. 
$$
This is a convex problem, so its minimum is computed by the Karush-Kuhn-Tucker 
theorem: there exists $\lambda \ge 0$ such that 
$\cL_{\bz}\, (\bz, \lambda) = 0$, where $\cL = - \sqrt{d}\prod_{i=1}^d z_i^{\,p_i} + \lambda 
\bigl( \sum_{i = 1}^d x_iz_i \, - \, 1\bigr)$. For each $i$, we have 
 $\cL_{\,z_i}\, = \, - \frac{p_i}{z_i}\sqrt{d}\prod_{i=1}^d z_i^{\,p_i} \, + \, 
 \lambda\, x_i \, = \, 0$. 
 Therefore $x_iz_i = c$ for all~$i$ . Substituting to the 
 constraint $\sum_{i=1}^d x_iz_i = 1$, we get $c = \frac{1}{d}$,  hence the point of minimum is 
 $z_i = 1/(d\,x_i)$. Consequently,  
 $$
 f^*(\bx) \ = \ \frac{(\bx, \bz)}{\sqrt{d}\prod_{i=1}^d z_i^{\,p_i}}\ = \ 
 \frac{1}{\sqrt{d}\prod_{i=1}^d (1/d\,x_i)^{p_i}}\ = \ 
 \sqrt{d}\prod_{i=1}^d x_i^{\,p_i} \ = \ f(\bx)\, .
 $$
This completes the proof. 
\smallskip 

{\hfill $\Box$}
\medskip

\begin{center}
\textbf{4.1. Basic properties of self-dual antinorms}
\end{center}
\bigskip

The geometrical meaning of  self-duality is that the unit ball 
${G \, =\, \{\bx \in \re^d_+ \, | \ f(\bx) \ge 1\}}$ is autopolar:~$G^{*} = G$, 
where the antipolar $G^*$ is defined in~(\ref{eq.antipolar}).
According to that definition, the polar to a point~$\bx$ is the half-space 
$\{\by \in \re^d \ | \ (\bx, \by) \ge 1\}$.  The  {\em polar hyperplane}~$p(\bx)$ is defined as the boundary of that subspace:  
 $p(\bx) = 
\{\by \in \re^d \ | \ (\bx, \by) = 1\}$. The point~$\bx$ is called the 
{\em pole} of this hyperplane.  We usually call the polar hyperplane 
simply {\em polar}, when it is clear that we mean the  plane but not the half-space. 
Let $\cS = \{\bx \in \re^d \ | \ f(\bx) = 1\}$ be the antisphere of $f$. 
The antisphere is {\it autopolar} if 
the set of its hyperplanes of support coincides with
$\{p(\bx) \, \bigl| \ \bx \in \cS\}$. This is equivalent to the autopolarity of the 
ball~$G$. The proof of the following fact 
is simple and we omit it. 
\begin{prop}\label{p.40}
An antinorm is self-dual if and only if its antisphere is autopolar. 
\end{prop}  

Another property of self-dual antinorms is that 
all of them  are smaller than the Euclidean norm in all but one direction, 
i.e., for every antinorm, there exists a 
unique direction where it is equal to the Euclidean norm. 
\begin{prop}\label{p.50}
For every self-dual antinorm~$f$, we have $f(\bx) \le |\bx|\, , \, \bx \in \re^d_+$, and there is a unique up to normalization vector $\ba \in \re^d_+, \, \ba \ne 0$, 
such that $f(\ba) = |\ba|$.   
\end{prop}  
{\tt Proof}. On the unit antisphere~$\cS$, we chose the point $\ba$
closest to the origin.
We have $\, 1 = f(\ba) = f^*(\ba) = \min_{f(\bx) > 0}\frac{(\bx, \ba)}{f(\bx)} \le  \frac{(\ba, \ba)}{f(\ba)} = (\ba, \ba)$. Hence $|\ba| \ge 1$. On the other hand, 
   if $|\ba| > 1$, then 
  the polar of the point $P$ strictly separates the whole antisphere $\cS$ from the origin. This contradicts to the self-duality. 
   Hence, $|\ba| = 1$ and $f(\bx) \le (\ba, \bx)$.
 Such a point~$\ba \in \cS$ is unique. Indeed, 
if there are two of them $\ba_1, \ba_2$, then for the point $\bc = \frac12(\ba_1 +\ba_2)$, 
we have $|\bc| < 1$ and $f(\bc) \ge 1$. Hence, the point $\frac{1}{f(\bc)} \bc$
belongs to $\cS$ and has the norm less than one, which is impossible.  

{\hfill $\Box$}
\medskip 

Thus, the vector $\ba$ in Proposition~\ref{p.50} is the closest point of the antispere $\cS$ to the origin and $|\ba| = 1$.
\begin{cor}\label{c.25}
For every self-dual antinorm, its antisphere has a unique 
point of intersection with the Euclidean unit sphere. 
This is the point~$\ba$ form Proposition~\ref{p.50}.  
\end{cor}

\bigskip

\begin{center}
\textbf{4.2. Classification of self-dual antinorms}
\end{center}
\bigskip

Once there are infinitely many self-dual antinorms in~$\re^d$, the question arises about 
their possible classification. For the plane ($d=2$) 
this problems admits a complete solution. First, we consider 
 the following
\smallskip 

\noindent \textbf{{\em Construction 1}}.
\smallskip 

\noindent 1) take arbitrary $\ba \in \re^2_+, \, |\ba| = 1$, and draw a ray 
$\{t\,\ba \ | \ t \in \re_+\}$, which  splits the positive orthant 
$\re^2_+$ into two 
angles $K_1$ and $K_2$ (one of them may be degenerate); 
\smallskip 

\noindent  2) choose an arbitrary antinorm~$f_1$ on $K_1$
such that $f_1(\ba) = 1$ and $f_1(\bx) \le (\ba,\bx), \, \bx \in K_1$; 
\smallskip 

\noindent 3) define the function $f_2$ on $K_2$ as a dual to $f_1$: 
$f_2(\bx_2) = \inf_{\bx_1 \in K_1}\frac{(\bx_1, \bx_2)}{f_1(\bx_1)}\, $.
Then define $f: \re^2_+ \to \re_+$ as  follows:  
\begin{equation}\label{eq.f}
f(\bx) \ = \ 
\left\{
\begin{array}{l}
f_1(\bx)\, , \, \bx \in K_1,\\ 
f_2(\bx)\, , \, \bx \in K_2\, . 
\end{array}
\right. 
\end{equation}

\begin{theorem}\label{th.20}
All self-dual antinorms in~$\re^2_+$ are precisely those 
defined by formula~(\ref{eq.f}) with functions $f_1, f_2$
obtained by Construction 1.   
\end{theorem} 
Before giving a proof we make one comment. 
\begin{remark}\label{r.22}
{\em 
The 
self-duality of the antinorm~$f$ defined in~(\ref{eq.f}) is not quite obvious, because 
$f_2$ is the dual to~$f_1$ only on the cone~$K_1$, but not in the entire~$\re^2_+$. 
To prove that actually~$f_2^* = f_1$ and vice versa one needs to show that 
for every~$\bx_2 \in K_2$, the minimum 
$\min_{\bx_1 \in \re^2_+}\frac{(\bx_1, \bx_2)}{f_1(\bx_1)}$ is attained 
for~$\bx_1 \in K_1$. }
\end{remark}

\begin{figure}[ht!]
\center{\includegraphics[width=0.35\linewidth]{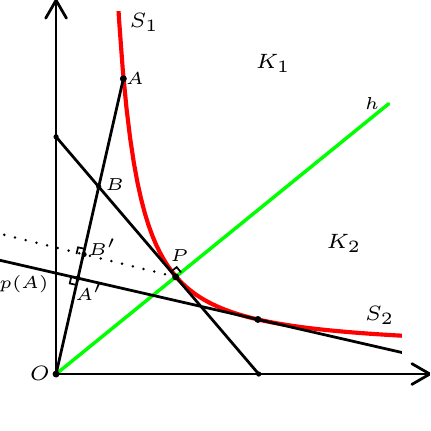}}
\caption{{\footnotesize Proof of Theorem~\ref{th.20}.}}
\label{pic1}
\end{figure}

{\tt Proof of Theorem~\ref{th.20}}. Let  $\cS_1 = \{\bx \in K_1 \ |\ f_1(\bx) = 1\}$
be the antisphere of the antinorm $f_1$ in $K_1$. Let $OP = \ba$. 
By the assumption, $P \in \cS_1$. For an arbitrary point $A \in \cS_1$, 
denote by $B$ the intersection of the segment $OA$ with the line orthogonal to 
$OP$ passing through $P$, and by $A', B'$ points inverse to $A, B$ respectively
(Fig.~\ref{pic1}). 
Since $|OA| \ge |OB|$ we have $|OA'| \le |OB'|$. Since in the right  triangle
$OPB$, the square of the leg $OP$ is equal to its projection to the 
hypotenuse $OB$ multiplied by $OB$, we see that the length of the projection is 
$|OP|^2/|OB| = 1/|OB| = |OB'|$. Hence $PB'$ is the polar $p(B)$. 
Therefore, the polar $p(A)$, which is parallel to $p(B)$ and passes through $A'$, 
is closer to $O$ than $p(B)$. 
Hence, the line $p(B)$ separates 
$p(A)$ from $\cS_1$. 
Thus, the polar $p(A)$ of an arbitrary point 
$A \in \cS_1$ does not intersect the line $\cS_1$ expect possibly at $P$. Hence, the whole polar image of $\cS_1$ is located in $K_2$ and therefore coincides with $\cS_2$. Then, by the reflexivity of the polar transform, 
the polar image of $\cS_1$ is located in $K_2$ and therefore coincides with $\cS_2$.
Thus, the polar transform interchanges $\cS_1$ and $\cS_2$, and $\cS$ is self-polar. 
   
   Conversely, let $f$ be a self-dual antinorm and $\cS$ be its antisphere. 
   By Proposition~\ref{p.50},  $\cS$ has a unique closest point $\ba$ to the origin and $|\ba| = 1$. 
   Denote $P = \ba$.  
    Let $h$ be the line orthogonal to $OP$ and passing through $P$. 
   The whole antisphere $\cS$ lies on the opposite side from the origin about $h$.
   Then the ray 
   $\{t\,\ba \, , t \in \re_+\}$ splits the orthant $\re^2_+$ into two angles 
   $K_1$ and $K_2$. Then the function $f_1 = f|_{K_1}$  generates the function 
   $f$ by the procedure 1-3 in the theorem.   
 \smallskip 

{\hfill $\Box$}
\medskip 

\begin{remark}\label{r.24}
{\em Geometrically, Theorem~\ref{th.20} gives the following characterization 
of self-dual antinorms on the plane. We split the positive orthant $\re^2_+$
by some ray~$\{t\boldsymbol{\alpha} \ | \ t \ge 0\}$ into two angles $K_1, K_2$ according to item 1) of Construction~1 
 and then take an arbitrary antinorm $f_1$
 in $K_1$ whose antisphere $\cS_1$ lies above the perpendicular  to that ray  
 through  $\ba$. Then $\cS_2$ is a polar of $\cS_1$. Concatenation of 
 $\cS_1$ and of $\cS_2$ gives the antisphere of $f$. Thus, the curve 
 $\cS_1$ defines the antinorm $f$. 
} 
\end{remark} 

 \begin{ex}\label{ex.18}
 {\em Let $\ba \, = \, \frac{\sqrt{2}}{2}\bigl(1, 1\bigr)$; then 
 the ray generated by~$\ba$ is the bisector of the coordinate angle~$XOY$. On this
 bisector we take a point $(R, R)$, where 
 $R = 1+\sqrt{2}$, and draw a circle of radius~$R$ centered at this point. 
 This circle is tangent to the axes and passes through~$\ba$. 
 We take as~$\cS_1$ the arc of this circle connecting the point~$\ba$ with the 
 point of tangency with the axis~$OY$. Then the polar~$\cS_2$ of this arc is 
 the piece of hyperbola $y \, = \, \frac{1}{x - a} + b$, 
 starting at~$\ba$ and going along the axis~$OX$ to $+\infty$, 
 where 
$a = \frac{3\sqrt{2}}{2} - 2, \, b =  \sqrt{2}-1$. It asymptotically tends to 
the line $y = \sqrt{2}-1$. Then the union of $\cS_1$ and $\cS_2$ is the 
autopolar antisphere.  
   }
 \end{ex}
 \smallskip 
 
 If $\cS_1$ is a broken line, then 
 $\cS_2$ is a broken line as well and $f$ is a piecewise linear function. 
 Hence, there 
 are  infinitely many   autopolar conic  polygons in~$\re^2_+$
 (in the sense of ``antipolar'') and we have classified them all. 
 In the next subsection we give their explicit description.

 \bigskip 

\newpage

\begin{center}
\textbf{4.3. Autopolar conic  polygons and polyhedra}
\end{center}
\bigskip 

The class of piecewise-linear antinorms is especially important, since such antinorms 
are easily described. A piecewise-linear  antinorm has the form 
 $f(\bx) \, = \, \min\limits_{j=1, \ldots, n} (\ba_j, \bx)$, where $\ba_j$ are 
nonnegative vectors. The unit ball of a piecewise-linear  antinorm is 
a conic polyhedron (see Definition~\ref{d.15}).  A {\em conic polygon} 
is a two-dimensional conic polyhedron . 

Theorem~\ref{th.20} allows us to provide  an explicit construction of 
all autopolar conic polygons. 

\smallskip 

\noindent \textbf{{\em Construction 2}}. 
We build a broken line $A_{-k}\ldots A_{-1}A_0A_1 \ldots A_{k-1}$  
in~$\re^2_+$ as follows. 
\medskip 

\begin{figure}[ht!]
\center{\includegraphics[width=0.35\linewidth]{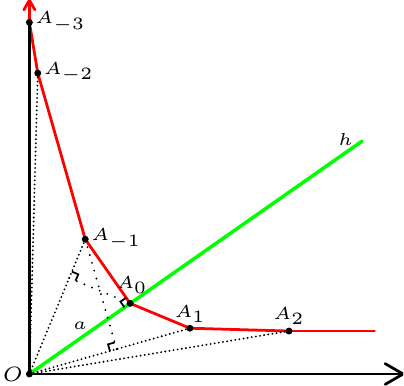}}
\caption{{\footnotesize Construction of an autopolar polygon for $k = 3$}}
\label{pic2}
\end{figure}

 Take an arbitrary point $A_0$ such that $OA_0 = 1$. 
Draw a line through $A_0$ orthogonal to $OA_0$ and take an arbitrary point 
$A_{-1}$ on it. Drop  a perpendicular from  $A_0$ to $OA_{-1}$ and take an arbitrary point 
$A_{1}$ on its extension  through~$A_0$. The following is
by induction. If the point $A_j$ is constructed (assume $j > 0$), then  we drop  a perpendicular from  $A_{-j}$ to $OA_j$ and take an arbitrary point~$A_{-j-1}$ 
 on its extension  through~$A_{-j}$. If $j< 0$, then we drop  a perpendicular from  $A_{-j-1}$ to $OA_j$ and take an arbitrary point~$A_{-j}$ 
 on its extension  through~$A_{-j-1}$ (Fig.\ref{pic2}). 

In the last iteration we choose the point $A_{-k}$ on the coordinate axis~$OY$. 
Draw a ray from $A_{-k}$ along the axis~$OY$ and a ray from 
$A_{k-1}$ parallel to the axis~$OX$;  
call those rays $A_{-k}Y$ and $A_{k-1}X$ respectively. 
Then $YA_{-k}\ldots A_{-1}A_0A_1 \ldots A_{k-1}X$ is the desired 
autopolar conic polygon.  
Invoking Theorem~\ref{th.20} we conclude

 \begin{cor}\label{c.28}
Every self-dual conic polygon is obtained by Construction 2.  
\end{cor} 
For each $k$, that conic polygon has $2k$ vertices (the vertex $A_{-k}$ is on the 
axis, all others are strictly inside~$\re^2_+$) and $2k+1$ sides 
(the sides $A_{-k}Y$ and  $A_{k-1}Y$ are  rays, all others are segments). 
The simplest cases are described in the following examples: 
 
 \begin{ex}\label{ex.20} ($k=0$). 
 {\em This case is formally out of our construction. Here  $A_0 = (1,0)$
 and the side $A_0A_{-1}$ becomes the vertical ray from~$A_0$. Thus, 
 the polygon has one vertex~$A_0$ and two sides that are rays from $A_0$
 parallel to the coordinate axes. The antinorm is $f(x,y) = x$. 
 }
 \end{ex}

 \begin{ex}\label{ex.25} ($k=1$). 
 {\em Here $A_0$ is an interior point of $\re^2_+$ and 
 $A_{-1}$ is a point on the axis $OY$ such that 
 $\angle \, A_{-1}A_0O \, = \, 90^{\circ}$. Then the self-dual 
 polygon is bounded by the rays $A_{-1}Y, A_0X$ and by the segment~$A_{-1}A_0$. 
 If $A_0 = (a,b)$, then $f(x,y) \, = \, \min \, \bigl\{ ax + by\, , \, \frac{y}{b}  \}$. 
  
 }
 \end{ex}
 
 Characterisation of self-dual antinorms, in particular, 
 polyhedral antinorms, in higher dimensions,  are left as open problems and discussed in Section~5.

\bigskip

\begin{center}
\textbf{4.4. Symmetric antinorms} 
\end{center}

\bigskip 
An antinorm is called {\em symmetric}
 if it is invariant with respect to every permutation of 
coordinates. There is a variety of symmetric antinorms. 
Choosing a parameter $p \in [-\infty, 1]$, we define 
the $L_p$-symmetrization of an arbitrary antinorm~$f$ as follows: 
$$
 f^{[p]}(x_1, \ldots , x_d) \ = \ \Bigl[\, \frac{1}{d!}\, \sum_{\sigma} f^{\,p}(x_{\sigma(1)} , \ldots , 
x_{\sigma(d)})\, \Bigr]^{1/p}, 
$$ with natural modifications for~$p=0$
(where the $L_p$-mean becomes the geometrical mean) and for $p=-\infty$
(the $L_p$-mean becomes the minimum). 
The sum  is computed over all 
permutations of the set~$\{1, \ldots , d\}$. Then~$f^{[p]}$ is a symmetric antinorm.

Are there symmetric self-dual antinorms? Yes, for example, 
$f(\bx) \ = \ \sqrt{d}\, \bigl(\, x_1\cdots x_d \, \bigr)^{1/d}$. 
This is the antinorm from Proposition~\ref{p.30} with all 
$p_i$ equal to~$1/d$.  What about others, do they exist?   In case of negative answer we come to the following situation: 
although there are many self-dual antinorms, in the class of symmetric antinorms 
it is unique. We are going to see that at least for $d=2$ this is true.  
 \smallskip

A bivariate  antinorm $f(x,y)$ on $\re^2_+$ is symmetric if 
$f(x, y) = f(y,x)$. We prove that in the rich variety of  self-dual antinorms 
on $\re^2_+$ (see Theorem~\ref{th.20}), there is only one symmetric antinorm. 
\begin{theorem}\label{th.30}
The unique symmetric self-dual antinorm on~$\re^2_+$ is  $f(x,y) = \sqrt{2xy}$. 
\end{theorem}
First, we realize the proof under the assumption that $f$ is smooth, i.e., 
is differentiable at each interior point of~$\re^2_+$.  
Geometrically  this means that  there is a unique line of support to the antiball $G = \{\bx \in \re^2_+ \ | \ 
f(\bx) \ge 1\}$ at every point of its boundary~$\cS = \partial B$. In this case we 
call the unique line of support {\em tangent line} to~$\cS$. Then we extend the proof to the nonsmooth case. 
We use two simple lemmas. 

\begin{lemma}\label{l.30}
If a monotone function $\varphi: \re_+ \to \re_+$ is such that 
$\varphi(0) = 0$ and $\varphi(\varphi(u))  \equiv u$, then  $\varphi(u)  \equiv u$.  
\end{lemma}
{\tt Proof.}  From the assumption it follows that $\varphi$ is non-decreasing. 
If $\varphi(u) >  u$ for some $u$, then taking $\varphi$ of both sides of the equation we obtain 
$\varphi(\varphi(u)) \ge \varphi(u) >  u$, hence 
$\varphi(\varphi(u))  >  u$, which is impossible. Similarly, if  $\varphi(u) <  u$, then 
$\varphi(\varphi(u)) <  u$. The contradiction proves that $\varphi(u) =  u$ for all $u \ge 0$.

{\hfill $\Box$}

\smallskip 

Now we are able to prove the following generalization of Lemma~\ref{l.30}:  
\begin{lemma}\label{l.40}
Let a set $\Omega \subset \re_+$ consist of disjoint nonempty half-open intervals
$(a_i, b_i], \, i \in \cI$, where the index set $\cI$ is either finite or coincides with $\n$.   
For each $u \in \re_+$, we define the number $\tilde u$ as follows: 
\begin{equation}\label{eq.con}
\tilde u  \ = \ 
\left\{
\begin{array}{lcl}
 u & , & u \notin \Omega \, , \\
 a_i & , & u \in (a_i, b_i], \ i \in \cI\ .  
 \end{array}
 \right.  
 \end{equation}
Then, if a monotone function ${\varphi: \re_+ \to \re_+}$ which is constant on each interval 
$(a_i, b_i], \, i \in \cI$, and  
$\varphi(0) = 0$,  satisfies the inequality $\varphi(\varphi(u)) \, = \, \tilde u$, 
then    $\varphi(u) \, = \, \tilde u$.

\end{lemma}
{\tt Proof.}  From the assumptions it follows  that 
 $\varphi$ is non-decreasing. If $\varphi(u) >   u$ for some~$u$, then 
 $\varphi(\varphi(u)) \ge   \varphi(u) > u$, which is 
 impossible, because~(\ref{eq.con}) implies that $\varphi(\varphi(u)) \le  u$
 for all~$u$. Consequently, $\varphi(u) \le   u$ for all~$u$. 
If $u \notin \Omega$, then $\varphi(\varphi(u)) = u$ and we conclude 
 as in the proof of Lemma~\ref{l.30} that  $\varphi(u) = u$.  If   $u \in (a_i, b_i]$, then 
$\varphi (u) \le u \le b_i$. If $\varphi (u) < a_i$, then 
$\varphi(\varphi (u)) \le \varphi(a_i) \le \varphi(u) <  a_i$, which is impossible due to~(\ref{eq.con}). 
Therefore, $\varphi(u) \in (a_i, b_i]$. Since $\varphi$
is constant on the interval $(a_i, b_i]$ and both $u$ and $\varphi(u)$
belong to it, we see that $\varphi(u) = \varphi(\varphi(u)) = a_i$, which completes the proof.

{\hfill $\Box$}

\begin{figure}[ht!]
\begin{minipage}[h]{0.32\linewidth}
\center{\includegraphics[width=1\linewidth]{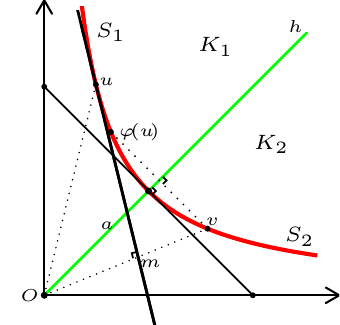}} a) \\
\end{minipage}
\hfill
\begin{minipage}[h]{0.32\linewidth}
\center{\includegraphics[width=1\linewidth]{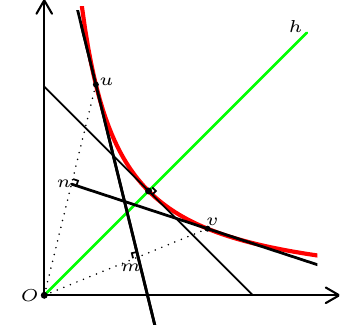}} b) \\
\end{minipage}
\hfill
\begin{minipage}[h]{0.32\linewidth}
\center{\includegraphics[width=0.97\linewidth]{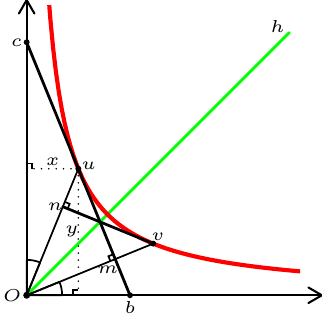}} c) \\
\end{minipage}
\caption{{\footnotesize Proof of Theorem~\ref{th.30}. The smooth case.}}
\label{pic5}
\end{figure}

{\tt Proof of Theorem~\ref{th.30}.} 
Let $f$ be a self-dual antinorm with unit antisphere $\cS$. By Proposition~\ref{p.50}, $|\bx| \ge 2$
for all $\bx \in \cS$, and there exists a unique $\ba \in \cS$ such that $|\ba| = 1$. 
If $f$ is
symmetric, then $\ba$ is symmetric to itself, and therefore $\ba = 
\bigl(\frac{1}{\sqrt{2}}, \frac{1}{\sqrt{2}} \bigr)$  (Fig.~\ref{pic5}\,a). 
\smallskip 

{\tt The smooth case.}
Define a map~$\varphi: \cS \to \cS$ as follows: 
for every $\bu \in \cS$, the point $\varphi(\bu)$ is symmetric to the pole 
of the tangent line to $\cS$ at the point~$\bu$. Denote this pole by~$\bv$ (Fig.~\ref{pic5}\,a).
Clearly, $\bv\in \cS$. The  self-duality implies that $\bu$ is a pole to the tangent line to~$\cS$ at the point $\bv$.  
Hence, $\varphi(\bv)$ is a point symmetric to~$\bu$. 

Let $\bu', \bv'$ be points symmetric to 
$\bu$ and~$\bv$ respectively (fig. 6 b). We have $\varphi(\bv) = \bu'$ and due to the symmetry, 
$\varphi(\bv') = \bu$. But $\bv' = \varphi(\bu)$. Therefore, 
$\varphi(\varphi (\bu)) = \bu$. The point $\ba$ splits the 
curve $\cS$ into two parts $\cS_1$ and $\cS_2$. 
Consider a natural parametrization of the
curve $\cS_1$ by the length of the arc from a point to $\ba$. 
It defines a homeomorphim of $\cS_1$ and $\re_+$. 
Let a point $\bu$ be associated to a parameter~$u$ and we keep 
the notation $\varphi$ for the corresponding map on the half-line. 
We have $\varphi (\varphi(u)) = u$ for all $u \in \re_+ $ and by 
Lemma~\ref{l.30}, $\varphi (u) = u$. Therefore, $\varphi (\bu) = \bu$, i.e., 
every point $\bu$ is symmetric to the pole~$\bv$ of the tangent line passing through~$\bu$.

  Denote by $\bm, \bb$, and $\bc$  the points of intersection of that tangent line 
  with the segment~$O\bv$ and with the coordinate axes $OX$ and $OY$ respectively. 
  We have $O\bm \cdot O\bv = 1$. Due to the symmetry, 
  $\angle \bb O \bv \, = \, \angle \bc O \bu$. On the other hand, 
  $O\bm$ is the altitude of the right triangle $\bb O \bc$ to the 
hypotenuse $\bb \bc$ (Fig.~\ref{pic5}\,c). Hence $\angle \bb O \bv \, = \, \angle  O \bc \bu$. 
Thus, $\angle  \bu O \bc \, = \, \angle  O \bc \bu$ and so $O\bu$ is 
the midpoint of the hypotenuse $\bb \bc$. Hence the median $O\bu$ is equal to the 
half of the hypotenuse $\bb \bc$. The area of the 
triangle $\bb O \bc$ is equal to $\frac12\, O\bm \cdot \bb \bc \, = \,  
O\bm \cdot O\bu \, = \, O\bm \cdot O\bv = 1$. On the other hand, the same 
area is equal to $\frac12 \, O\bb \cdot O\bc \, = \, \frac12 \, 2x\cdot 2y \, = \, 
 2xy$, where $x,y$ are the coordinates of the point $\bu$ (Fig.~\ref{pic5}\,c). Thus, $2xy = 1$
 for every point $\bu \in \cS$. Hence,  $f(x,y) = 1$ if and only if $2xy = 1$, 
 and by homogeneity  $f(x,y) = \sqrt{2xy}$.

\begin{figure}[ht!]
\center{\includegraphics[width=0.4\linewidth]{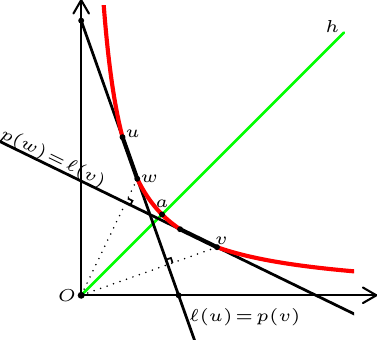}}
\caption{{\footnotesize Proof of Theorem~\ref{th.30}. The non-smooth case. }}
\label{pic6}
\end{figure}

 {\tt The non-smooth case.}
Define a map~$\varphi: \cS \to \cS$ as follows: 
for every $\bu \in \cS$, 
draw a line of support to $\cS$ at the point~$\bu$
closest to the point~$\ba$. Denote this straight line by~$\ell(\bu)$. 
Then define $\varphi(\bu)$  
as the point symmetric (with respect to the bisector) to the pole~$\bv$ 
of the line $\ell(\bu)$ (Fig.~\ref{pic6}).

 Clearly, if $\cS$ has an edge, then  all 
points $\bu$ of that edge, except for its end closest to~$\ba$, have 
the same image $\varphi(\bu)$. Indeed, for all such points~$\bu$,  
the line~$\ell(\bu)$ is the same. Therefore it suffices to consider the case 
when $\bu$ is the end of the edge most distant from~$\ba$. 
  Denote by $\bw$ the other end of that edge.
If $\bw = \bu$, then the edge~$[\bu, \bw]$ is trivial and  $\cS$ is smooth at the point~$\bu$. Since $f$ is self-dual, we have $\bv \in \cS$ and 
the set of  polars to all 
$\bx \in [\bu, \bw]$ is the set of lines of support to $\cS$ at the point~$\bv$. 
Among them, the polar to $\bw$ is the closest one to $\ba$. 
Hence, the pole of the line $\ell(\bv)$ is  $\bw$ and consequently 
 $\varphi(\bv) = \bw'$, where  $\bw', \bu', \bv'$ are the points symmetric 
to $\bw, \bu$, and $\bv$ respectively about the bisector of the coordinate angle. 
Since $\cS$ is symmetric, it contains all $\bw', \bu'$, and $\bv'$ 
and $\varphi(\bv') = \bw$. Taking into account that 
$\bv' = \varphi(\bu)$, we conclude that $\varphi\bigl(\varphi(\bu)\bigr) \, = \, \bw$.
In particular, if $f$ is smooth at $\bu$, then $\varphi\bigl(\varphi(\bu)\bigr) \, = \, \bu$. 

Now consider the natural parametrisation of $\cS_1$ and 
the corresponding function $\varphi: \re_+ \to \re_+$. It satisfies all 
the assumptions of Lemma~\ref{l.40} with the set $\Omega$ being 
the image of the  union of nontrivial half-open edges of~$\cS_1$
(we remove the end closest to~$\ba$ from each edge). 
 Applying Lemma~\ref{l.40} we 
obtain $\varphi(\bu) \, = \, \bw$. Hence $\bw$ is symmetric to $\bv$. 
Consider the right triangle $\bb O\bc$ formed by the line $\ell(\bu)$
and by the two coordinate axes. 
As in the proof for the smooth case we conclude that 
$\bw$ is the midpoint  of the hypotenuse $\bb \bc$ and  that the area of $\bb O\bc$  is $2xy$, where $\bw = (x,y)$
(Fig.~\ref{pic6}).
This is well known that in this case 
the line $\ba \bb$ touches  the hyperbola $\cH = \{(x, y) \in \re^2_+ \ | \ 2xy = 1\}$ at the point~$\bw$. Actually, $\cH$ is not the whole hyperbola, 
but only one of two branches,  but we 
keep a short notation.

 Thus,  we have proved the following property of the curve $\cS$: 
{\em each edge  of this curve is tangent to the hyperbola $\cH$ 
at its end closest to the point~$\ba$}    (Fig.~\ref{pic6}). In particular, 
all points of $\cS$ out of nontrivial edges belong to~$\cS$. This implies that 
 $\cS = \cH$. Indeed, if $\cS$ has a proper edge $[\bu, \bw]$, 
then it touches $\cH$ at the point $\bw$ and hence $\bu \notin \cH$. 
Therefore, there exists a point $\bx \in \cS$ close to $\bu$ such  that 
the arc of $\cS$ between the points $\bu$ and $\bx$ intersects 
 neither  $\cH$ nor the half-open interval $(\bu , \bw]$. Since  $\bx$ does not lie on~$\cH$, it 
 must belong to a proper edge~$[\bm, \bn]$. In this case $\bn \in \cH$, which is impossible, 
 since $\bn$ belongs to the arc between $\bu$ and $\bx$, which does not intersect~$\cH$.     
The contradiction proves that $\cS$ does not have proper edges, therefore all its 
points belong to $\cH$ and so $\cS = \cH$.

{\hfill $\Box$}

\bigskip 

Theorem~\ref{th.30} implies that the only autopolar symmetric conic body in 
$\re^2_+$ is the hyperbola~$\sqrt{2xy} = 1$. 
Possible generalizations of this result to $\re^d_+$
are discussed in the next section.

\bigskip 

\begin{center}
\large{\textbf{5. Open problems}}
\end{center}

\bigskip 

As shown in Section~4, there are infinitely many self-dual antinorms in 
$\re^d_+$. However, we succeed in their  classifying 
only for $d=2$. So, the first open problem is the following: 
\bigskip

\noindent \textbf{Problem~1}. {\em How to characterise 
 self-dual antinorms in~$\re^d_+$ for $d\ge 3$? }
\bigskip 

In the two-dimensional case, every self-dual antinorm is constructed  starting with the 
vector~$\ba$ (see Proposition~\ref{p.50}). Draw the ray  $h\,  = \, \{t\ba \ | \ t \ge 0\}$ and choose arbitrary mutually dual antinorms in those parts. 
In the $d$-dimensional case, Proposition~\ref{p.50} still holds and we can draw a ray $h$, but it does not split the positive orthant~$\re^d_+$. Most likely, constructing 
self-dual antinorms for~$d\ge 3$ requires other ideas.  
\bigskip 

The situation with polyhedral self-dual antinorms is still more complicated. 
As for now, we are not only able to classify them for~$d\ge 3$ but 
do not know if they exist at all, apart from  liftings of two-dimensional antinorms. 
  The lifting is defined as follows. 
Let $P$ be a $k$-dimensional 
conic polyhedron 
in some of $k$-dimensional faces of $\re^d_+$, say, in the face 
consisting of points  $(x_1, \ldots , x_k, 0, \ldots , 0)$. 
The {\em lifting} of $P$ in $\re^d$  is 
 $G = \bigl\{(x_1, \ldots , x_d) \in \re^d_+\ | \ (x_1, \ldots , x_k) \in P
\bigr\}$. Clearly, $G$ is also a conic polyhedron, which is the 
right cylinder with the base~$P$. Moreover, if $P$ is autopolar, then so 
is~$G$. If $P$ is a unit ball 
of the antinorm $\varphi(x_1, \ldots , x_k)$ on the corresponding face of~$\re^d_+$, then 
   $G$ is a unit ball 
of the antinorm $f (\bx) = \varphi(x_1, \ldots , x_x)$ on~$\re^d_+$. 
 
 Each autopolar conic polygon in~$\re^2_+$ produces an autopolar
 conic polyhedron  in~$\re^d_+$ by the lifting. Hence, 
 there are infinitely many $d$-dimensional conic polyhedra. 
 The question is if there are others? 
 \smallskip  

 \noindent \textbf{Problem~2}. {\em Does there exist at least one autopolar conic polyhedron in $\re^3$, 
 which is not a lifting of a smaller-dimensional conic polyhedron? The same question is for 
higher dimensions.}  
 \medskip 
 
In two-dimensional case we built autopolar conic polygons starting with the 
point~$\ba$ (the closest point to the origin) and then 
constructing successively other vertices, see subsection~4.3. In higher dimensions, however,  
this approach does not seem to be applicable.  
 
 \bigskip 
 
 Our last problem concerns the symmetric antinorms. 
  \bigskip 
  
 \noindent \textbf{Problem~3}. {\em Is it true that $f(\bx) \ = \ \sqrt{d}\, \bigl(\, x_1\cdots x_d \, \bigr)^{1/d}$ is the unique self-dual symmetric antinorm in 
 $\re^d_+$?}
\smallskip  
 
 For $d=2$, the affirmative answer is given in Theorem~\ref{th.30}.   
 Its proof was based on the general construction of two-dimensional 
 self-dual antinorms, which is not applicable in higher dimensions.  
 It is interesting whether the proof of the $d$-dimensional analogue 
 of  Theorem~\ref{th.30}, provided it is correct, can be derived 
 by applying Theorem~\ref{th.30} and a kind of inductive argument? 
 
\bigskip

\begin{center}
\large{\textbf{6. Applications}}
\end{center}
\bigskip 
 
We consider several applications of antinorms: the linear switching systems, 
the largest Lyapunov exponent, the lower spectral radius of matrices, and the convex trigonometry.
We will see that the role and the interpretation of duality and of self-duality is different in those applications. 
 \bigskip

 \begin{center}
\textbf{6.1. Lyapunov functions for linear switching systems}
\end{center}
\bigskip

The linear switching system is a dynamical system defined for an arbitrary 
compact set of $d\times d$ matrices~$\cA$ (the {\em control set}) as follows:   
\begin{equation}\label{eq.lls}
\left\{
\begin{array}{l}
\dot \bx(t) \ = \ A(t)\bx(t)\, , \qquad t \ge 0\\
\bx(0)\ = \ \bx_0\\
A(t) \ \in \ \cA \, , \qquad   t \ge 0
\end{array}
\right. 
\end{equation}
The function~$A(\cdot)$ called {\em the switching law} is a measurable 
function from $\re_+$ to the control set of matrices~$\cA$. Thus, we have 
a linear differential equation with the controlled matrix $A(\cdot)$ 
in the right hand side. 

See~\cite{L03, MP89} for general theory and applications
of linear switching systems.  
One of important issues  is estimating  
the fastest possible growth of trajectories~$\bx(t)$
as $t\to \infty$. A measure of the fastest growth~$\hat \sigma (\cA)$ is 
  equal to the infimum 
of numbers $\alpha$ such that $\|\bx(t)\| \, \le \, C \, e^{\, \alpha t}$
for every trajectory~$\bx(t)$.  The system is {\em (asymptotically) stable} 
if every trajectory 
tends to zero  as $t \to +\infty$.  The stability  is equivalent to the condition
 ${\hat \sigma (\cA) < 0}$~\cite{MP86}. The standard way to prove the stability is to 
 find a {\em Lyapunov function}, which is a continuous positive homogeneous 
 function on~$\re^d$ that decreases along every trajectory. The existence of a Lyapunov function implies the stability. 
 The converse is also true. Moreover, each stable system possesses a convex
 Lyapunov function, i.e., a norm decreasing along every trajectory~\cite{MP86}.
 The most popular Lyapunov functions (quadratic, piecewise-quadratic, polyhedral, 
 sum-of-squares, etc.) are norms.  If the control set~$\cA$ is convex, 
  then there exists an {\em invariant norm} also called~{\em Barabanov norm}~$\|\cdot \|$ such that for every 
  trajectory, we have 
  $\|\bx(t)\| \le e^{\hat \sigma} \|\bx(0)\|, \, t \ge 0$, 
  and for every point~$\bx_0$, there exists a trajectory 
   $\bar \bx(t)$ starting at that point and such that 
  $\|\bx(t)\| = e^{\hat \sigma} \|\bx(0)\|, \, t \ge 0$~\cite{B88}.  
  Apart from some special cases, it is not known how to find invariant norms. 
  Nevertheless, constructing ``nearly invariant''  norms makes it possible to 
  compute~$\hat \sigma$   with a good precision~\cite{BS08, GLP17,  LA09, MP89}.

Similarly, the slowest growth of trajectories of linear switching systems is characterized by 
the  value  $\check \sigma (\cA)$ equal to  the infimum
of numbers $\alpha$, for which there exists a control function $A(\cdot)$ such that
for every starting point~$\bx_0 \ne 0$,  the corresponding trajectory satisfies
$\|\bx(t)\| \, \le \, C \, e^{\, \alpha t}$.  The system is {\em stabilizable} if
there is a control function $A(\cdot)$ such that every trajectory 
with the switching law~$A(\cdot)$
tends to zero  as $t \to +\infty$ independently of the initial point~$\bx_0$.  The stabilizability is equivalent to the condition
 $\check \sigma < 0$~\cite{LA09, SDP08}. To decide the stabilizability one can 
 consider a homogeneous positive Lyapunov function~$f: \re^d \to \re_+$
 which increases along   every trajectory. If such a function exists, then clearly the system is 
 not stabilizable. However, the converse is in general not true: for 
 non stabilizable systems,  such a function  may not exist~\cite{BS08}.  The reason is that 
 here we cannot rely on convexity. All known proofs of the existence of the 
 Lyapunov norm for stable systems~\cite{B88, W02}  use the fact that the pointwise supremum of convex
 functions is also  convex. For stabilizable systems, one needs to use infinum instead of supremum, but for infimum this property does not hold: the pointwise infimum does not respect convexity.  
 One could then  replace the convexity by concavity. 
 However, there are no concave positive functions on~$\re^d$ apart from identical constants. 
 Nevertheless, such functions exist on cones. In fact,
 if the system has an invariant cone, then there always exists 
 a Lyapunov concave function (antinorm) on that cone. This is true, in particular, for {\em positive systems} that have their trajectories in the positive orthant~$\re^d_+$.
Moreover, a positive system always  possesses an {\em invariant antinorm}~$f$ such that for every 
  trajectory, we have 
  $f(\bx(t)) \ge e^{\check \sigma} f(\bx(0)), \, t \ge 0$, 
  and for every point~$\bx_0 \in \re^d_+$, there exists a trajectory 
   $\bar \bx(t)$ starting at that point and such that 
  $f(\bar \bx(t)) = e^{\check \sigma} f(\bar \bx(0))$ for all $\, t \ge 0\, $~\cite{GLP17}.
  Stability and stabilizabilty of positive systems have been studied in many works~\cite{FM12, GLP17} and references therein. 
  
  Approximating the invariant antinorm makes it possible to 
  compute the Lyapunov exponent~$\hat \sigma$   with a prescribed precision. 
  Duality of antinorms enables us to find the following relations between 
  the systems with control sets~$\cA$ and $\cA^T = \{A^T \, | \, A \in \cA\}$. 
  \begin{prop}\label{p.60}
  If $f$ is a Lyapunov antinorm for the system with a control set~$\cA$, then 
  $f^*$ is a Lyapunov antinorm for the system with the control set~$\cA^T$. 
  \end{prop}
  \begin{remark}\label{r.30}
{\em Proposition~\ref{p.60} is not that obvious 
as it seems. For example, it is not true for invariant antinorms: 
if $f$ is an invariant antinorm for~$\cA$, then $f^*$ is, in general, not an  invariant antinorm 
for~$\cA^T$
}
\end{remark}
{\tt Proof of Proposition~\ref{p.60}}. If $f$ is a Lyapunov antinorm, then for every sufficiently small 
 $s > 0$,  there exists $\varepsilon  = \varepsilon(s) > 0$ such that 
$f(\bx +  s A\bx) \, < \, (1 - \varepsilon)\, f(\bx)$ for every $\bx > 0$
and $A \in \cA$~\cite{GLP17}. Then, for every $\by > 0$, we have 
$$
f^*\bigl(\by - s A^T\by \bigr) \ = \ \inf_{\bx > 0} \frac{(\bx , (I - s A^T)\by)}{f(\bx)} \ = \ 
\inf_{\bx > 0} \frac{((I - s A)\bx , \by)}{f(\bx)} \ \le \ 
\inf_{\bx > 0} \frac{((I - s A)\bx , \by)}{f(\bx)}\, . 
$$
Denote $\bz = (I - s A)\bx$. Then $\bx = (I - s A)^{-1}\bz \, = \, (I + s A + o(s))\bz$
as $s \to 0$. Hence, for sufficiently small~$S$, we have 
$$
f^*\bigl(\by - s A^T\by \bigr) \ = \ 
\inf_{\bz > 0} \frac{(\bz , \by)}{f( \bz +  s A\bz + o(s))}\ >  
\ \inf_{\bz > 0} \frac{(\bz , \by)}{(1 - \varepsilon)f(\bz)}\ = \ \frac{1}{1 - \varepsilon} f^*(\by)\ > \ (1 + \varepsilon)f^*(\by).   
$$
Thus, $f^*\bigl(\by - s A^T\by \bigr) \, > \,  (1 + \varepsilon)f^*(\by)$
for every~$A^T \in \cA^T$, 
therefore~$f^*$ is a Lyapunov antinorm for the family~$\cA^T$, see~\cite{GLP17}.

{\hfill $\Box$}
\medskip 

Proposition~\ref{p.60} establishes the duality of Lyapunov antinorms 
for transposed families. In practice, this allows us to 
construct the  Lyapunov function by passing to the transpose family, 
for which this problem is sometimes simpler.

\bigskip

 \begin{center}
\textbf{6.2. The largest Lyapunov exponent of matrices}
\end{center}
\bigskip 

The antinorm is also applied for computation of the largest Lyapunov exponent
arising in the multiplicative ergodic theorem. For the sake of simplicity, we consider the case of 
products of independent random matrices distributed over a finite matrix family. 
Let us have a family of $d\times d$-matrices $\cA = \{A_1, \ldots , A_m\}$. 
To each matrix~$A_j$ we associate a positive probability~$p_j$
so that $\sum_{j=1}^m p_j= 1$.  Consider a random product  $X_k = A_{d_k}\cdots A_{d_1}$,
where all indices $\{d_j\}$ are independent and identically distributed  random variables; each $d_j$
takes values $1, \ldots , m$ with probabilities $p_1, \ldots , p_m$ respectively. According to
the Furstenberg-Kesten theorem~\cite{FK60} the value $\|X_k\|^{1/k}$ converges with probability~$1$
to a number $r$, which depends only on the family~$\cA$,  and on the probabilities
$\{p_j\}_{j=1}^m$. The number $\ell = \log r$ is called the {\em largest Lyapunov exponent} of~$\cA$. A strong generalization of this theorem was proved by Oseledets~\cite{O68}. 
We do not deal with other Lyapunov exponents, and omit the word ``largest''.  This, however, should not lead  
to a confusion with the Lyapunov exponent for the linear switching system 
(from the previous section). This number can be found by the following
 limit formula
 \begin{equation}\label{LE}
 \ell  \quad = \quad \lim_{k \to \infty} \ \frac1k \ \bE \log \ \bigl\| \, A_{d_k}\cdots A_{d_1} \bigr\|\, ,
 \end{equation}
where $\bE$ denotes the mathematical expectation. The computation of  the Lyapunov exponent is hard even for $2\times 2$ matrices. No efficient algorithms are known. This is not surprising because the Lyapunov exponent is in general a discontinuous function of matrices. 
Moreover, the problem to distinguish between two cases: $\ell \ge 0$ and $\ell <0$ is 
algorithmically undecidable~\cite{BT97}. For nonnegative matrices, the situation is 
slightly better. In  this case the Lyapunov exponent is continuous, and there are  numerically efficient  
algorithm for its approximate  computation, 
see, for instance~\cite{H97,  JP13, P08} and references therein.   
 Although the  distinguishing between 
 $\ell \ge 0$ and $\ell <0$ is NP-hard even for Boolean matrices~\cite{BT97}.
  For properties of Lyapunov exponents of non-negative matrices, 
  see~\cite{H97, P13, W86}.

 Some of the computational algorithms are based on approximation of the 
Lyapunov antinorm~$f$ on~$\re^d_+$, which is characterised by the property: 
$$
\prod_{j=1}^m f^{p_j}(A_j\bx)  \ <  \  f(\bx)\, , \qquad \bx \in \re^d_{+}\, . 
$$
Under some mild assumption on matrices (the absence of zero rows and columns and of common invariant coordinate subspaces), the following holds:  $\ell < 0$ if and only if there exists a strictly positive Lyapunov antinorm~\cite{P10}.   Moreover, under the same assumptions, there always exists 
an {\em invariant antinorm}, for which 
$$
\prod_{j=1}^m f^{p_j}(A_j\bx)\ = \ r\, f(\bx) \,  , \qquad \bx \in \re^d_{+}\, , 
$$
see~\cite{P11}. 
Surprisingly enough, an analogue of Proposition~\ref{p.60} is not true for the 
random matrix products: the dual to  a Lyapunov antinorm 
of a matrix family~$\cA$ may not be a Lyapunov antinorm 
for the transpose family~$\cA^T$. 
 \begin{prop}\label{p.70}
There exists a family~$\cA$ of nonnegative $2\times 2$ matrices 
  and an antinorm $f$ on~$\re^2_+$ such that $f$ is a Lyapunov 
 function for~$\cA$, but $f^*$ is not a  Lyapunov 
 function for~$\cA^T$. 
  \end{prop}
{\tt Proof}. Take an arbitrary $q \in \bigl(\frac{\sqrt{2}}{2}, 1 \bigr)$ and consider the 
following pair of matrices: 
$$
A_1 \ = \ 
q\, \left(
\begin{array}{cc}
1 & 1 \\
0 & 1
\end{array}
\right)\, ; \qquad 
A_2 \ = \ 
q\, \left(
\begin{array}{cc}
1 & 0 \\
1 & 1
\end{array}
\right)\, 
$$
with the probabilities $p_1 = p_2 = \frac12$. Since $A_2 = A_1^T$, the family~$\cA = \{A_1, A_2\}$ satisfies $\cA = \cA^T$. 
Consider the antinorm~$f(x, y) = x+y$ on~$\re^2_+$. Then 
its dual is $f^*(x, y) = \min\, \{x, y\}$ (see Example~\ref{ex.12}). 
For each point  $\bx = (x, y)$ such that $f(x, y) = 1$, we have 
$A_1\bx \, = \, q(x+y, y)^T\, = \, q(1, y)^T$. Hence, $f(A_1\bx)\, = \, q(1+y)$. 
Similarly, $A_2\bx \, = \, q(x, 1)^T$ and $f(A_2\bx)\, = \, q(1+x)$. 
Therefore, 
$$ 
f^{p_1}(A_1\bx)f^{p_2}(A_2\bx)\ = \ 
\sqrt{f(A_1\bx)f(A_2\bx)} \ = \ q\sqrt{(1+y)(1+x)} \ \ge \ q \sqrt{2}  \ > \ 
1 \ = \ f(\bx)
$$
Thus, for every $\bx \in \re^2_+$, we have $f^{p_1}(A_1\bx)f^{p_2}(A_2\bx) \, > \, 
f(\bx)$, 
hence $f$ is a Luapunov antinorm for~$\cA$. 

On the other hand, if $f*(\bx) = \min\, \{x, y\} = 1$, 
then $f^*(A_1\bx) \, = \, \min\, \{q(x+y),  qy\}\, = \, qy$
and $f^*(A_2\bx) \, = \, \min\, \{qx,  q(x+y)\}\, = \, qx$, consequently 
$\sqrt{f^*(A_1\bx)f^*(A_2\bx)} \, = \, q \sqrt{xy}$. At the point~$\bx = (1, 1)$, 
we have $\sqrt{f^*(A_1\bx)f^*(A_2\bx)} \, = \, q \, < \, 1$, and so 
$\sqrt{f^*(A_1\bx)f^*(A_2\bx)} \,  < \, f^*(\bx)$, therefore, 
$f^*$ is not a Lyapunov antinorm for~$\cA$, neither for $\cA^T = \cA$.

{\hfill $\Box$}
\medskip 

Thus, dual antinorms may not correspond to 
the Lyapunov antinorms of transposed families of matrices. 
It would be interesting to understand the sense of duality for 
the   Lyapunov antinorms.

\bigskip

 \begin{center}
\textbf{6.3. The lower spectral radius}
\end{center}
\bigskip 

\smallskip

The {\em lower spectral radius} (also called in the literature the {\em joint spectral 
subradius})  of a compact family of  matrices~$\cA$ is 
$$
\check \rho(\cA) \ = \ \lim_{k\to \infty} \min_{d_1, \ldots , d_k}\|A_{d_k}\cdots A_{d_1}\|^{1/k}\, , 
$$
where the minimum is defined over all possible products of length~$k$ of matrices from 
$\cA$, with repetitions permitted. This limit always exists and does not depend on the 
matrix norm.  For a family of one matrix, the lower spectral radius becomes the 
usual spectral radius of that matrix, which is the maximum modulus of its eigenvalues. 
The lower spectral radius is the exponent of the minimal growth 
of matrix products of length $k$ as $k\to \infty$. 

The lower spectral radius was introduced in~\cite{G95} to characterise the 
minimal growth of trajectories  of the discrete-time 
linear switching system~$\bx(k+1) = A(k)\bx(k)$ as $k\to \infty$, where 
$A(k) \in \cA, \, k \ge 0$. The sequence~$A(k)$ is called the {\em switching law}. 
The system is {\em stabilizable} if there 
exists a switching law whose trajectory~$\{\bx(k)\}_{k \ge 0}$ tends to zero 
as $k \to \infty$, for every initial point~$\bx(0) \in \re^d$.
The stabilizability is equivalent to the condition~$\check \rho (\cA) < 1$. 
See~\cite{BM14, GP13, M17} for more properties of the lower spectral radius. 
Apart from the dynamical systems, it has found applications in 
the theory of wavelets, in approximation theory, in the number theory, combinatorics, the 
theory of formal languages, etc.  

Many of those applications 
(see, for example,~\cite{FV12, FV13, GP13, JPB09, P17}) deal with 
nonnegative matrix families.
For them, the lower spectral radius can be efficiently bounded in terms of antinorms. 
For an arbitrary nonnegative family~$\cA$, there exists an {\em extremal}
antinorm~$f$ on~$\re^d_+$ such that
$$
\min_{A \in \cA}\, f(A\bx) \ \ge  \ \check \rho (\cA) \, f(\bx)\, \quad \bx \in \re^d_+\, . 
$$
Moreover, under some mild assumptions on~$\cA$, there is also an {\em invariant antinorm}~\cite[Theorems 5,6]{GP13}: 
\begin{defi}\label{d.50} 
An antinorm~$f$ on~$\re^d_+$ is called invariant for a family~$\cA$ 
of non-negative $d\times d$ matrices if 
$$
\min_{A \in \cA}\, f(A\bx) \ =  \ \check \rho(\cA) \, f(\bx)\, \quad \bx \in \re^d_+\, .  
$$
\end{defi}
Efficient algorithms for approximate computation of the lower spectral radius 
are based on the construction of  extremal and invariant 
antinorms~\cite{GP13}. Those  antinorms are also of an independent interest~\cite{GZ15}.

The following theorem establishes a relation between the  extremal antinorms of a matrix family and 
of its transpose. Thus, for the lower spectral radius, this relation is found, 
unlike for the Lyapunov exponent of random matrix products (Proposition~\ref{p.70}). 
 Moreover, it is also possible to characterize not only the extremal antinorm, but also the invariant antinorm of a transpose family, in contrast to the situation with the 
 continuous-time switching systems~(Remark~\ref{r.30}). To formulate the theorem we 
 need to introduce some more notation. For an arbitrary subset $X$ of $\re^d_+$, 
 we consider its {\em positive convex hull}~$\, {\rm co}_+\, X\, = \, 
 {\rm co}\, X\,  + \, \re^d_+\, = \, 
 \bigl\{\bx \in \re^d_+\ | \ \exists \, \by \in {\rm co}\, X\, \bigr\}$, 
 where, recall, ${\rm co}\, X$ is the (standard) convex hull of~$X$.   
 A conic body $P \subset \re^d_+$ is called {\em invariant}
 for a family~$\cA$ of nonnegative matrices if 
 ${\rm co}_+\, \bigl\{\, AP \ | \ A \in \cA\, \bigr\}\,  =  \, \check \rho \, P$. 
 \begin{theorem}\label{th.40}
If an antinorm~$f$ is extremal for $\cA$, then 
$f^*$ is extremal for~$\cA^T$. If an antinorm~$f$ is invariant for $\cA$, 
then the polar $G^*$ to its unit ball is an invariant conic body for~$\cA^T$, 
and vice versa: if $G$ is an invariant conic body for~$\cA$, then 
the antinorm with the unit ball~$G^*$ is invariant for~$\cA^T$. 
  \end{theorem}
{\tt Proof}.  Without
loss of generality it can be assumed that $\check \rho(\cA) = 1$. 
If~$f$ is an extremal antinorm for $\cA$, then 
\begin{equation}\label{eq.100}
\min_{\bx \ge 0, \, f(\bx) \ge 1\, , A \in \cA}  \ f(A\bx)\quad = \quad 1\, . 
\end{equation}
For the dual antinorm $f^*$, we have 
$\ \min\limits_{\bp \ge 0, \, \, f^*(\bp) \ge 1\, , A \in \cA}  \ f^*(A^T\bp)\ = $
\begin{equation}\label{eq.110}
\min_{\bp \ge 0, \, \, f^*(\bp) \ge 1\, , A \in \cA} \quad \min_{\bx \ge 0, f(\bx) \ge  1}\ (A^T\bp, \bx)\ = \ 
\min_{\bx, \bp \ge 0, \, \, f(\bx) \ge 1, \,  f^*(\bp) \ge 1\, , A \in \cA} \quad  (\bp, A\bx) \, . 
\end{equation}
Applying~(\ref{eq.100}), and replacing  $\by = A\bx$, we see that the last expression in~(\ref{eq.110}) is equal to 
$$
\min_{\bx, \bp \ge 0, \, \, f(A\bx) \ge 1, \,  f^*(\bp) \ge 1\, , A \in \cA} \quad (\bp, A\bx)\ =\ 
\min_{\by, \bp \ge 0, \, \, f(\by) \ge 1, \,  f^*(\bp) \ge 1} (\bp, \by)\ = \ 
$$
$$
\min_{\bp \ge 0,  f^*(\bp) \ge 1}\quad \min_{\by \ge 0,  f^*(\by) \ge 1} \ (\bp, \by) \quad  = 
\quad  
\min_{\bp \ge 0,  f^*(\bp) \ge 1}\  f^*(\bp)\ = \ 1.  
$$
Thus, $\min\limits_{\bp \ge 0, \, f^*(\bp) \ge 1\, , A^T \in \cA^T}  \ f^*(A^T\bp)\ = \ 1$, 
hence, $f^*$ is an extremal antinorm for $\cA^T$, which completes the proof of the first 
statement. 

Now we need to prove that  if $f$ is an invariant antinorm for $\cA$, 
then $G^*$  in an invariant conic body for~$\cA^T$, 
and vice versa. We prove only the primal assertion, after which the dual assertion becomes 
obvious. If $f$ is an invariant antinorm for $\cA$, then, 
by definition, for every~$\bx \in \re^d_+$,  
we have $\min_{A\in \cA} f(A\bx)\, = \, f(\bx)$. 
Moreover, the function of support 
of the set $A^TG^*$ at the point~$\bx$ is equal to $f(A\bx)$. Indeed, 
$$
\min_{\by \in A^TG^*} \, (\bx, \by) \ = \ \min_{\bp \in G^*} \, (\bx, A^T\bp) \ = \
 \min_{\bp \in G^*} \, (A\bx, \bp)\ = \ f^{**}(A\bx)\ = \ f(A\bx), 
$$
where the latter equality follows from Corollary~\ref{c.20} and from continuity 
of~$f$. Similarly, the function of support 
of $G^*$ at the point~$\bx$ is equal to $f(\bx)$. Since
$\min_{A\in \cA} f(A\bx)\, = \, f(\bx)$, we see that the minimal 
function of support of the sets~$A^TG^*$ over all~$A \in \cA$ is equal to 
the function of support of~$G^*$. Hence, 
${\rm co}_+\, \{A^TG^* \ | \ A\in \cA\}\, = \, G^*$. 

{\hfill $\Box$}
\medskip 

\begin{remark}\label{r.40}
{\em It is interesting that while the duality takes an extremal antinorm of matrices to  
an extremal antinorm of their transposes, it does not  do the same with the invariant antinorm. 
If $f$ is an invariant antinorm for~$\cA$, then $f^*$ is not invariant 
for~$\cA^T$ but $G^*$  is an invariant conic body for~$\cA^T$. For the joint spectral radius, which is responsible for the stability of discrete-time linear switching system, 
a similar relation between the invariant norm and the invariant body of the transpose family 
was established in~\cite{PW08}.   
}
\end{remark}

\bigskip 

\begin{center}
\textbf{6.4. The convex trigonometry}
\end{center}
\bigskip 

One more application of the duality and self-duality of antinorms 
is in the extension of convex trigonometry to the convex hyperbolic functions. 
The convex trigonometry was developed recently~\cite{L19} in the study of bivariate optimal control problems. Let us recall the basic construction. For an arbitrary convex body~$G \subset \re^2$ 
containing the origin~$O$ in its interior, we denote by $\cS$ its boundary, by 
$\mathbb{S}$ its area, and for an arbitrary number $\theta \in [0, 2\mathbb{S}]$, we 
define by~$P_{\theta}$ the point on~$\cS$ such that the oriented area 
of the part of~$G$ bounded by rays $OX$ and $OP_{\theta}$ is equal to~$\frac12 \theta$. 
This definition is extended to all~$\theta \in \re$ in a standard manner by periodicity. 
Then by definition $\cos_{G}\, \theta $ and  $\sin_{G}\, \theta $ are respectively the 
abscissa and the ordinate of $P_{\theta}$. The standard trigonometric functions 
correspond to the case when $G$ is a unit disc. It was shown~\cite{L19} that all basic 
trigonometric formulas can be generalized to $\cos_{G}\, \theta $ and  $\sin_{G}\, \theta $ . 
In particular, the identity $\cos^{2}\theta \, + \, \sin^{2}\theta \, = \, 1 $, 
becomes $\cos_{G}\theta\, \cos_{G^*}\theta^* \, + \, \sin_{G}\theta\, \sin_{G^*}\theta^* \, = \, 1 $, where $G^*$ is the polar of~$G$ and $\theta^*$ corresponds to 
the point on the boundary of~$G^*$ defined by the direction of the normal 
to~$G$ drawn at the point~$P_{\theta}$. 

Thus, every convex body~$G \subset \re^2$ containing the origin as an interior point defines trigonometric functions. In particular, every norm in~$\re^2$ defines them 
by means of its unit ball. It was announced in~\cite{L19}  that a similar construction 
can produce hyperbolic functions. Indeed, they can be defined by an arbitrary 
antinorm in~$\re^2_+$ (in this case it should be rather called 
{\em concave trigonometry}). To see this, we consider an antinorm~$f$ and 
the closest to the origin  point~$P_0$ of  the unit sphere~$\cS$ of~$f$. 
Denote by  $\tilde \cS$ the unit sphere~$\cS$ 
complemented by rays of the coordinate axes
in case $\cS$ intersects the corresponding axes. 
Thus, if $f$ vanishes on the axes, then $\tilde \cS = \cS$. 

 For an arbitrary $\theta \in \re$, 
we consider the point~$P_{\theta} \in \tilde \cS$
such that the oriented area 
of the part of the unit ball~$G$ bounded by $\tilde \cS$, by the line~$OP_0$, 
and by the perpendicular dropped from~$P_{\theta}$ to 
that line, is equal to~$\frac12 \theta$. 
By definition $\cosh_{G}\, \theta $ and  $\sinh_{G}\, \theta $ are respectively the 
abscissa and the ordinate of $P_{\theta}$. The standard hyperbolic functions 
correspond to the case when $G$ is the hyperbola $2xy = 1$. 
The identity $\cosh^{2}\theta \, - \, \sinh^{2}\theta \, = \, 1 $, 
becomes $\cosh_{G}\theta\, \cosh_{G^*}\theta^* \, + \, \sinh_{G}\theta\, \sinh_{G^*}\theta^* \, = \, 1 $, where $G^*$ is the polar of~$G$ and $\theta^*$ corresponds to 
the point on the boundary of~$G^*$ defined by the direction parallel to the normal 
to~$G$ drawn at the point~$P_{\theta}$. If~$G$ is autopolar, i.e., $f$
is self-dual, then this formula is simplified to
$\cosh_{G}\theta\, \cosh_{G}\theta^* \, + \, \sinh_{G}\theta\, \sinh_{G}\theta^* \, = \, 1 $. 
In contrast to the convex trigonometry, here we have a variety of 
autopolar sets classified by~Theorem~\ref{th.30}. 
On the other hand, only one of them is symmetric. This is the hyperbola $2xy =1$
corresponding to the antinorm $f(x, y) = \sqrt{2xy}$ and producing the standard 
hyperbolic functions. 
\vspace{1cm}

\noindent \textbf{Acknowledgements.}  The author is grateful to the anonymous referee for his impressive work and for many valuable comments. He also expresses his  thanks to N.Guglielmi for 
many useful discussions and to T.Zaitseva for her help in illustrations.

 \end{document}